\documentclass[preprint,3p,sort&compress,final,times]{elsarticle}

\usepackage{amssymb}

\usepackage{placeins}
\usepackage{tikz}
\usepackage{ams math,cases}
\usepackage{siunitx}
\usepackage{mathtools}

\usepackage{subcaption}
\usepackage{placeins}

\newcommand{\subsubfloat}[2]{%
  \begin{tabular}{@{}c@{}}#1\\#2\end{tabular}%
}

\newcommand{\figref}[1]{Figure~\ref{#1}} 
\newcommand{\secref}[1]{Section~\ref{#1}} 

\definecolor{darkred}{rgb}{0.5,0,0}
\definecolor{darkgreen}{rgb}{0,0.3,0}
\definecolor{darkblue}{rgb}{0,0,0.5}
\definecolor{darkbrown}{rgb}{0.28,0.07,0.07}
\definecolor{black}{rgb}{0,0,0}
\definecolor{plotGreen}{rgb}{0,0.5,0}
\definecolor{plotRed}{rgb}{1.0,0,0}
\definecolor{plotBlue}{rgb}{0,0,1.0}
\definecolor{plotCyan}{rgb}{0,1.0,1.0}
\definecolor{plotGray}{gray}{0.25}

\newcommand{\velocity}{\vec{u}}
\newcommand{\velocityd}{\vec{u}_{h}}
\newcommand{\velocitydv}{\boldsymbol{u}}
\newcommand{\velocitydc}[1]{u_{#1}}
\newcommand{\velocitydnext}{\velocityd^{\,k+1}}
\newcommand{\velocitydprevious}{\velocityd^{\,k}}
\newcommand{\velocitydmida}{\tilde{\vec{u}}_{h}^{\,k+\frac{1}{2}}}
\newcommand{\velocitydmidb}{\vec{u}_{h}^{\,k+\frac{1}{2}}}
\newcommand{\velocitydmidnext}{\vec{u}_{h}^{\,k+\frac{3}{2}}}

\newcommand{\velocitydvmidb}{\boldsymbol{u}^{\,k+\frac{1}{2}}}
\newcommand{\velocitydvmidnext}{\boldsymbol{u}_{h}^{\,k+\frac{3}{2}}}
\newcommand{\velocitydmidstart}{\vec{u}_{h}^{\,\frac{1}{2}}}
\newcommand{\velocitydstart}{\vec{u}_{h}^{\,0}}
\newcommand{\velocitydstartnext}{\vec{u}_{h}^{\,1}}
\newcommand{\velocitydreverseprevious}{\vec{u}_{h}^{\,k+\frac{1}{2}}}
\newcommand{\velocitydreversenext}{\vec{u}_{h}^{\,k+\frac{3}{2}}}

\newcommand{\vorticity}{\omega}
\newcommand{\vorticityd}{\omega_{h}}
\newcommand{\vorticitydv}{\boldsymbol{\omega}}
\newcommand{\vorticitydc}[1]{\omega_{#1}}
\newcommand{\vorticitydnext}{\vorticityd^{\,k+1}}
\newcommand{\vorticitydprevious}{\vorticityd^{\,k}}
\newcommand{\vorticitydvnext}{\vorticitydv^{\,k+1}}
\newcommand{\vorticitydvprevious}{\vorticitydv^{\,k}}
\newcommand{\vorticitydmida}{\tilde{\omega}_{h}^{\,k+\frac{1}{2}}}

\newcommand{\vorticitydstart}{\vorticityd^{\,0}}
\newcommand{\vorticitydreverseprevious}{\vorticityd^{\,k}}
\newcommand{\vorticitydreversenext}{\vorticityd^{\,k+1}}

\newcommand{\gentranspd}{\vorticityd}
\newcommand{\gentranspdnext}{\vorticitydnext}
\newcommand{\gentranspdprevious}{\vorticitydprevious}

\newcommand{\ufunction}{\vec{v}_{h}}
\newcommand{\pfunction}{q_{h}}
\newcommand{\wfunction}{\xi_{h}}
\newcommand{\ufunctionc}{\vec{v}}
\newcommand{\pfunctionc}{q}
\newcommand{\wfunctionc}{\xi}
\newcommand{\uspace}{H(\mathrm{div},\Omega)}
\newcommand{\udspace}{U_{h}}
\newcommand{\ubasis}[1]{\vec{\epsilon}^{\,U}_{#1}}
\newcommand{\udim}{d_{U}}
\newcommand{\wspace}{H(\mathrm{curl},\Omega)}
\newcommand{\wdspace}{W_{h}}
\newcommand{\wbasis}[1]{\epsilon^{W}_{#1}}
\newcommand{\wdim}{d_{W}}
\newcommand{\pspace}{L^{2}(\Omega)}
\newcommand{\pdspace}{Q_{h}}
\newcommand{\pbasis}[1]{\epsilon^{Q}_{#1}}
\newcommand{\pdim}{d_{Q}}
\newcommand{\totalp}{\bar{p}}
\newcommand{\totalpd}{\bar{p}_{h}}
\newcommand{\totalpdv}{\boldsymbol{\bar{p}}}
\newcommand{\totalpdc}[1]{p_{#1}}
\newcommand{\totalpdnext}{\totalpd^{\,k+1}}

\newcommand{\totalpdprevious}{\totalpd^{\,k}}
\newcommand{\totalpdvnext}{\totalpdv^{\,k+1}}

\newcommand{\totalpdmida}{\tilde{\bar{p}}_{h}^{\,k+\frac{1}{2}}}

\newcommand{\rtspace}{\mathrm{RT}_{N}}
\newcommand{\qspace}{\mathrm{DG}_{N-1}}
\newcommand{\cgspace}{\mathrm{CG}_{N}}

\newcommand{\energy}{\mathcal{K}}
\newcommand{\energyd}{\mathcal{K}_{h}}
\newcommand{\energydnext}{\mathcal{K}_{h}^{k+\frac{3}{2}}}
\newcommand{\energydprevious}{\mathcal{K}_{h}^{k+\frac{1}{2}}}

\newcommand{\enstrophy}{\mathcal{E}}
\newcommand{\enstrophyd}{\mathcal{E}_{h}}
\newcommand{\enstrophydnext}{\mathcal{E}_{h}^{k+1}}
\newcommand{\enstrophydprevious}{\mathcal{E}_{h}^{k}}

\newcommand{\vorticityt}{\mathcal{W}}
\newcommand{\vorticitytd}{\mathcal{W}_{h}}
\newcommand{\vorticitytdnext}{\mathcal{W}_{h}^{k+1}}
\newcommand{\vorticitytdprevious}{\mathcal{W}_{h}^{k}}

\newcommand{\matrixoperator}[1]{\boldsymbol{\mathsf{#1}}}
\newcommand{\matrixoperatorc}[1]{\mathsf{#1}}

\begin{document}

\begin{frontmatter}



\journal{Journal of Computational Physics}
\title{A mass, energy, enstrophy and vorticity conserving (MEEVC) mimetic spectral element discretization for the 2D incompressible Navier-Stokes equations}


\author[cst]{A.~Palha\corref{cor}}
\ead{a.palha@tue.nl}
\author[tudelft]{M.~Gerritsma}

\cortext[cor]{Corresponding author}

\address[cst]{Eindhoven University of Technology, Department of Mechanical Engineering, P.O. Box 513, 5600 MB Eindhoven, The Netherlands}

\address[tudelft]{Delft University of Technology, Faculty of Aerospace Engineering, P.O. Box 5058, 2600 GB Delft, The Netherlands}

\begin{abstract}
	In this work we present a mimetic spectral element discretization for the 2D incompressible Navier-Stokes equations that in the limit of vanishing dissipation exactly preserves mass, kinetic energy, enstrophy and total vorticity on unstructured grids. The essential ingredients to achieve this are: (i) a velocity-vorticity formulation in rotational form, (ii) a sequence of function spaces capable of exactly satisfying the divergence free nature of the velocity field, and (iii) a conserving time integrator. Proofs for the exact discrete conservation properties are presented together with numerical test cases on highly irregular grids.

\end{abstract}


\end{frontmatter}

\section{Introduction}
	\subsection{Relevance of structure preserving methods}
		Structure-preserving discretizations are known for their robustness and accuracy. They conserve fundamental properties of the equations (mass, momentum, kinetic energy, etc). For example, it is well known that the application of conventional discretization techniques to inviscid flows generates artificial energy dissipation that pollutes the energy spectrum. For these reasons, structure-preserving discretizations have recently gained popularity. 
		
		For a long time, in the development of general circulation models used in weather forecast, it has been noticed that care must be taken in the construction of the discretization of physical laws. Phillips \cite{Phillips1959} verified that the long-time integration of non-linear convection terms resulted in the breakdown of numerical simulations independently of the time step, due to the amplification of weak instabilities. Later, Arakawa \cite{Arakawa1966} proved that such instabilities can be avoided if the integral of the square of the advected quantity is conserved (kinetic energy, enstrophy in 2D, for example). Staggered finite difference discretizations that avoid these instabilities have been introduced both by Harlow and Welch \cite{Harlow1965} and Arakawa and his collaborators \cite{Arakawa1977, Mesinger1976}. Lilly \cite{Lilly1965} showed that these discretizations could conserve momentum, energy and circulation. At the same time, Piacseck and Williams \cite{Piacsek1970} connected the conservation of energy with the preservation of skew-symmetry of the convection operator at the discrete level. These ideas of staggering the discrete physical quantities and of preserving the skew-symmetry of the convection operator have been successfully explored by several authors aiming to construct more robust and accurate numerical methods. 
		
		In this work we focus on the incompressible Navier-Stokes equations, particularly convection-dominated flow problems (e.g. wind turbine wake aerodynamics).  It is widely known that in the absence of external forces and viscosity these equations contain important symmetries or invariants, e.g. \cite{Arnold1966,Arnold1992,Majda2001,Foias2001}, such as conservation of kinetic energy. A straightforward discretization using standard methods does not guarantee the conservation of these invariants. Discrete energy conservation is important from a physical point of view, especially for turbulent flow simulations when Direct Numerical Simulation (DNS) or Large Eddy Simulation (LES) are used. In these cases, the accurate reproduction of the energy spectrum is essential to generate the associated energy cascade. The numerical diffusion inherent (and sometimes essential) to many discretizations can dominate both the molecular diffusion contribution (DNS) and the sub-grid model contribution (LES). This negatively affects the energy spectrum and consequently the energy cascade. Energy-conserving discretizations ensure that all diffusion is modelled and not a product of discretization errors. For this reason, many authors have shown that energy-conserving schemes are essential both for DNS (e.g. \cite{Perot1993,Verstappen1995,Le1997,Verstappen1998,Verstappen2003}) and LES simulations (e.g. \cite{Mahesh2004,Mittal1997,Benhamadouche2002,Nagarajan2003,Felten2006,Ham2007}). In addition, discrete kinetic energy conservation provides a non-linear stability bound to the solution (e.g. \cite{Sadourny1975,Sanderse2013}). This bound relaxes the characteristic stability constraints of standard methods, allowing the choice of mesh and time step size to be based only on accuracy requirements. This is particularly relevant for LES where mesh sizes have to be kept as large as possible, due to the still significant computational effort required for this approach. Energy-conserving methods generate well-behaved global errors and adequate physical behaviour, even on coarse meshes.
	
	In two-dimensions, enstrophy is another conserved quantity of the flow when external forces and viscosity are not present. The classical Arakawa scheme \cite{Arakawa1966} and derived schemes exploit both the conservation of energy and enstrophy. By doing so, these methods have far better long time simulation properties than competing methods that do not conserve these two quantities.
		
	\subsection{Overview of structure preserving methods}	
		Two of the most used approaches to construct energy preserving discretizations are: (i) staggering and (ii) skew-symmetrization of the convective operator.
		
		Most staggered grid methods date back to the pioneering work of Harlow and Welch \cite{Harlow1965}, and Arakawa and colleagues \cite{Arakawa1977, Mesinger1976}. These methods employ a discretization that distributes the different physical quantities (pressure, velocity, vorticity, etc) at different locations in the mesh (vertices, faces, cell centres). It can be shown that, by doing so, important conservation properties can be maintained. Due to its success, much attention has been given to this approach and several extensions to the original work have been made by several authors. Morinishi \cite{Morinishi1998} has developed a high-order version on Cartesian grids. Extensions to non-uniform meshes have been presented for example in \cite{Wesseling1999} for structured quadrilateral meshes and by several authors \cite{Apanovich1988,Choudhury1990,Perot2000,Mullen2009} for simplicial meshes.
		
		The advantages of exploiting the skew-symmetry of the convection operator at the discrete level, were first identified by Piacseck \cite{Piacsek1970}. Following his ideas, Verstappen and Veldman \cite{Verstappen,Verstappen1998,Verstappen2003} and later Knicker \cite{Knikker2009} constructed high order discretizations on Cartesian grids with substantially improved properties. Kok \cite{Kok2009} presents another formulation valid on curvilinear grids. The extension to simplicial grids has been reported for example in the works of Vasilyev \cite{Vasilyev2000}, Veldman \cite{Veldman2008} and van't Hof \cite{VanHof2012}.
		
		At another level, Perot \cite{perot43discrete} suggested that a connection exists between the skew-self-adjoint formulation presented by Tadmor \cite{Tadmor1984} and box schemes. Box schemes or Keller box schemes are a class of numerical methods originally introduced by Wendroff \cite{wendroff1960} for hyperbolic problems and later popularized by Keller \cite{keller1971,keller1978} for parabolic problems. This method is face-based and space and time are coupled by introducing a space-time control volume. This method is known to be physically accurate and successful application has been reported in \cite{Croisille2002,Croisille2005,Gustafsson2006,Ranjan2013}. More recently, this method has been shown to be multisymplectic by Ascher \cite{Ascher2005} and Frank \cite{Frank2006}. Perot, \cite{Perot2007}, also established a relation between box schemes and discrete calculus and generalized it to arbitrary meshes.
		
		Although most of the literature on the simulation of incompressible Navier-Stokes flows has been performed using finite difference and finite volume methods, finite element methods have also been actively used and shared many of the developments already discussed for finite differences and finite volumes. One of the initial challenges in using finite elements for flow simulations lies in the fact that specific finite element subspaces must be used to discretize the different physical quantities. The finite element subspaces must satisfy the Ladyzhenskaya-Babuska-Brezzi (LBB) (or inf-sup) compatibility condition (see \cite{brezzi1991mixed}), otherwise instability ensues. Several families of suitable finite elements have been proposed in the literature, the most common being the Taylor-Hood family (Taylor and Hood \cite{TaylorHood1973}) and the Crouzeix-Raviart family (Crouzeix and Raviart \cite{Crouzeix}). The underlying idea behind these different finite element families is to use different polynomial orders for velocity and pressure. In essence, this is intimately related to the staggering approach already discussed in the context of finite differences and finite volumes. Staggering is explicitly mentioned in some finite element work, for example \cite{KoprivaKolias, Liu2007,Liu2008,Chung2012,Tavelli2015}. Another important aspect when using finite elements is the weak formulation used for the convection term. Some of these forms (rotational and skew-symmetric) can retain important symmetries of the original system of equations, this has been reported for example in \cite{Guevremont1990,Blaisdell1996,Ronquist1996,Wilhelm2000}. Ensuring conservation in finite element discretizations is not straightforward, nevertheless examples in the literature exist, e.g. \cite{Liu2000,Bernsen2006}. Also within the context of finite elements, Rebholz and co-authors, e.g. \cite{Rebholz2007,Olshanskii2010}, have developed several velocity-vorticity discretizations for the 3D Navier-Stokes equations. These discretizations are capable of conserving both energy and helicity.
		
		One final aspect when discussing existing structure preserving methods is time integration. Many of the conserving methods reported in the literature present proofs regarding spatial discretization. Nevertheless, when discrete time evolution is taken into account many time integrators destroy the nice properties of the spatial discretization. The relevance of time integration for the construction of structure preserving discretizations has been analyzed by Sanderse \cite{Sanderse2013}.

	\subsection{Overview of mimetic discretizations}
		Over the years numerical analysts have developed numerical schemes which preserve some of the structure of the differential models they aim to approximate, so in that respect the whole structure preserving idea is not new. One of the contributions of mimetic methods is to identify the proper language in which to encode these structures/symmetries is the language of differential geometry. Another novel aspect of mimetic discretizations is the identification of the metric-free part of differential models, which can (and should) be conveniently described in terms of algebraic topology. A general introduction and overview to spatial and temporal mimetic/geometric methods can be found in \cite{Christiansen2011,perot43discrete,Budd2003,Hairer2006}.

		The relation between differential geometry and algebraic topology in physical theories was first established by Tonti \cite{tonti1975formal}. Around the same time Dodziuk \cite{Dodziuk76} set up a finite difference framework for harmonic functions based on Hodge theory. Both Tonti and Dodziuk introduce differential forms and cochain spaces as the building blocks for their theory. The relation between differential forms and cochains is established by the Whitney map ($k$-cochains $\rightarrow$ $k$-forms) and the de Rham map ($k$-forms $\rightarrow$ $k$-cochains). The interpolation of cochains to differential forms on a triangular grid was already established by Whitney, \cite{Whitney57}. These interpolatory forms are now known as the {\em Whitney forms}.

		Hyman and Scovel \cite {HymanScovel90} set up the discrete framework in terms of cochains, which are the natural building blocks of finite volume methods. Later Bochev and Hyman \cite{bochev2006principles} extended this work and derived discrete operators such as the discrete wedge product, the discrete codifferential, the discrete inner product, etc.

		In a finite difference/volume context Robidoux, Hyman, Steinberg and Shashkov, \cite{HymanShashkovSteinberg97,HymanShashkovSteinberg2002,HYmanSteinberg2004,RobidouxAdjointGradients1996,RobidouxThesis,bookShashkov,Steinberg1996,SteibergZingano2009} used symmetry considerations to discretize diffusion problems on rough grids and with non-smooth anisotropic diffusion coefficients. In a more recent paper by Robidoux and Steinberg \cite{Robidoux2011} a discrete vector calculus in a finite difference setting is presented. Here the differential operators grad, curl and div are exactly represented at the discrete level and the numerical approximations are all contained in the constitutive relations, which are already polluted by modeling and experimental error. For mimetic finite differences, see also Brezzi et al. \cite{BrezziBuffaLipnikov2009,brezzi2010}.

		The application of mimetic ideas to unstructured staggered grids has been extensively studied by Perot, \cite{Perot2000,ZhangSchmidtPerot2002,perot2006mimetic,PerotSubramanian2007a,PerotSubramanian2007}. Especially in \cite{perot43discrete} where the rationale of preserving symmetries in numerical algorithms is lucidly described. The most \emph{geometric approach} is described in the work by Desbrun et al. \cite{desbrun2005discrete,ElcottTongetal2007,MullenCraneetal2009,PavlovMullenetal2010} and the thesis by Hirani \cite{Hirani_phd_2003}. 
		
		The \emph{Japanese papers}  by Bossavit, \cite{bossavit_japan_computational_1,bossavit_japan_computational_2,bossavit_japan_computational_3,bossavit_japan_computational_4,bossavit_japan_computational_5}, serve as an excellent introduction and motivation for the use of differential forms in the description of physics and the use in numerical modeling. The field of application is electromagnetism, but these papers are sufficiently general to extend to other physical theories. 
		
		In a series of papers by Arnold, Falk and Winther, \cite{arnold:Quads,arnold2006finite,arnold2010finite}, a finite element exterior calculus framework is developed. Higher order methods are also described by Rapetti \cite{Rapetti2007,Rapetti2009} and Hiptmair \cite{hiptmair2001}. Possible extensions to spectral methods were described by Robidoux, \cite{robidoux-polynomial}. A different approach for constructing arbitrary order mimetic finite elements has been proposed by the authors \cite{Palha2014, gerritsmaicosahom2012,Rebelo2014,palhaAdvectionIcosahom2014,kreeft::stokes,bouman::icosahom2009,palha::icosahom2009}.
		
		Extensions of these ideas to polyhedral meshes have been proposed by Ern, Bonelle and co-authors in \cite{Bonelle2015,Bonelle2014,Bonelle2015a,Bonelle2016} and by Brezzi and co-authors in \cite{BeiraodaVeiga2014,Brezzi2014,BeiraodaVeiga2016,DaVeiga2015}. These approaches provide more geometrical flexibility while maintaining fundamental structure preserving properties.
		
	Mimetic isogeometric discretizations have been introduced by Buffa et al. \cite{BuffaDeFalcoSangalli2011}, Evans and Hughes \cite{Evans2013a}, and Hiemstra et al. \cite{Hiemstra2014}.
	
	Another approach develops a discretization of the physical field laws based on a discrete variational principle for the discrete Lagrangian action. This approach has been used in the past to construct variational integrators for Lagrangian systems, e.g. \cite{Kouranbaeva2000,Marsden2003}. Recently, Kraus and Maj \cite{Kraus2015} have used the method of formal Lagrangians to derive generalized Lagrangians for non-Lagrangian systems of equations. This allows to apply variational techniques to construct structure preserving discretizations on a much wider range of systems.

	\subsection{Outline of paper}
		In this work we present a new exact mass, energy, enstrophy and vorticity conserving (MEEVC) mimetic spectral element solver for the 2D incompressible Navier-Stokes equations. The essential ingredients to achieve this are: (i) a velocity-vorticity formulation in rotational form, (ii) a sequence of function spaces capable of exactly satisfying the divergence free nature of the velocity field, and (iii) a conserving time integrator. This results in a set of two decoupled equations: one for the evolution of velocity and another one for the evolution of vorticity.
		
		In \secref{sec::spatial_discretization} we present the spatial discretization. We start by introducting the $(\velocity,\vorticity)$ formulation based in the rotational form in \secref{subsec::the_v_omega_formulation_in_the_rotational_form} and then in \secref{subsec::finite_element_discretization} the finite element discretization is discussed. This is followed by the temporal discretization in \secref{sec::temporal_discretization}. Once we have introduced the numerical discretization its conservation properties are proved in \secref{sec::conservation_properties_and_time_reversibility}. In \secref{sec::numerical_test_cases} the method is applied and tested on two test cases. We start by testing the accuracy of the method with a Taylor-Green vortex, for which the analytical solution is known. We finalize the test cases with an inviscid shear layer roll-up test to illustrate the conservation properties of the method. In \secref{sec::conclusions} we conclude with a discussion of the merits and limitations of this solver and future extensions.
	
\section{Spatial discretization} \label{sec::spatial_discretization}
	
	\subsection{The $(\velocity,\vorticity)$ formulation in rotational form} \label{subsec::the_v_omega_formulation_in_the_rotational_form}
		The evolution of 2D viscous incompressible flows is governed by the Navier-Stokes equations, which are most commonly expressed as a set of conservation laws for momentum and mass involving the velocity $\velocity$ and pressure $p$:
		\begin{equation}
			\begin{dcases}
				\frac{\partial\velocity}{\partial t} + \left(\velocity\cdot\nabla\right)\velocity + \nabla p = \nu\Delta\velocity + \vec{s}, \\
				\nabla\cdot\velocity = 0,
			\end{dcases} \label{eq::ns_convective_form}
		\end{equation}
		with $\nu$ the kinematic viscosity, $\vec{s}$ the body force per unit mass, and $\Delta = \nabla\cdot\nabla$ the Laplace operator. These equations are valid on the fluid domain $\Omega$, together with suitable initial and boundary conditions. In this work we consider only periodic boundary conditions and we set $\vec{s}=0$.
		
		The form of the Navier-Stokes equations presented in \eqref{eq::ns_convective_form} is the so called \emph{convective form}. Its name stems from the particular form of the nonlinear term $ \left(\velocity\cdot\nabla\right)\velocity$, which underlines its convective nature. This form is not unique. Using well known vector calculus identities it is possible to rewrite the nonlinear term in three other forms, see for example Zang \cite{Zang1991}, Morinishi \cite{Morinishi1998}, and R{\o}nquist \cite{Ronquist1996}.
		
		The first alternative, called \emph{divergence form} or \emph{conservative form}, expresses the nonlinear term as a divergence:
		\begin{equation}
			\begin{dcases}
				\frac{\partial\velocity}{\partial t} + \nabla\cdot\left(\velocity\otimes\velocity\right)+ \nabla p = \nu\Delta\velocity, \\
				\nabla\cdot\velocity = 0.
			\end{dcases} \label{eq::ns_divergence_form}
		\end{equation}
		
		The second alternative is obtained as a linear combination of the convective and divergence forms, producing a skew-symmetric nonlinear term, thus its name \emph{skew-symmetric form}:
		\begin{equation}
			\begin{dcases}
				\frac{\partial\velocity}{\partial t} + \frac{1}{2} \left(\velocity\cdot\nabla\right)\velocity + \frac{1}{2}\nabla\cdot\left(\velocity\otimes\velocity\right)+ \nabla p = \nu\Delta\velocity, \\
				\nabla\cdot\velocity = 0.
			\end{dcases} \label{eq::ns_skew_symmetric_form}
		\end{equation}
		
		The third and final alternative we consider here, the \emph{rotational form}, makes use of the vorticity $\vorticity:=\nabla\times\velocity$:
		\begin{equation}
			\begin{dcases}
				\frac{\partial\velocity}{\partial t} + \omega\times\velocity + \nabla \totalp = \nu\Delta\velocity, \\
				\nabla\cdot\velocity = 0,
			\end{dcases} \label{eq::ns_rotational_form}
		\end{equation}
		where the static pressure $p$ has been replaced by the total pressure $\totalp := \frac{1}{2}\velocity\cdot\velocity + p$.
		
		A different, but equivalent, approach used to describe fluid flow problems resorts to expressing the momentum equation in terms of the vorticity $\vorticity$. By taking the curl of the momentum equation and using the kinematic definition $\vorticity := \nabla\times\velocity$ we can obtain the flow equations based on vorticity transport:
		\begin{equation}
			\begin{dcases}
				\frac{\partial\vorticity}{\partial t} + \frac{1}{2}\left(\velocity\cdot\nabla\right)\vorticity + \frac{1}{2}\nabla\cdot\left(\velocity\,\vorticity\right) = \nu\Delta\vorticity, \\
				\nabla\cdot\velocity = 0, \\
				\vorticity = \nabla\times\velocity\,.
			\end{dcases} \label{eq::ns_vorticity_transport}
		\end{equation}
		This velocity-vorticity $(\velocity,\vorticity)$ formulation of the Navier-Stokes equations is of particular interest for vortex dominated flows, see for example Gatski \cite{Gatski1991} for an overview and Daube \cite{Daube1992} and Clercx \cite{Clercx1997a} for applications.
		
		All these different ways of expressing the governing equations of fluid flow are equivalent at the continuous level. As stated before, one can start with \eqref{eq::ns_convective_form} and derive all other sets of equations simply by using well known vector calculus identities and definitions. The often overlooked aspect is that  each of these formulations leads to a different discretization, with its own set of properties, e.g. \cite{Zang1991,Gatski1991,Ronquist1996,Morinishi1998}. Only in the limit of vanishing mesh size ($h \rightarrow 0$) are the discretizations expected to be equivalent. This introduces an additional degree of freedom associated to the choice of formulation that, combined with the discretization approach to use, will determine the final properties of the method.
		
		In this work we construct a $(\velocity,\vorticity)$ formulation by combining the rotational form \eqref{eq::ns_rotational_form} with the vorticity transport equation \eqref{eq::ns_vorticity_transport}, similar to the work of Benzi et al. \cite{Benzi2012} and Lee et al. \cite{Lee2011}:
		\begin{equation}
			\begin{dcases}
				\frac{\partial\velocity}{\partial t} + \omega\times\velocity + \nabla \totalp = -\nu\nabla\times\vorticity, \\
				\frac{\partial\vorticity}{\partial t} + \frac{1}{2}\left(\velocity\cdot\nabla\right)\vorticity + \frac{1}{2}\nabla\cdot\left(\velocity\,\vorticity\right) = \nu\Delta\vorticity, \\
				\nabla\cdot\velocity = 0\,,
			\end{dcases} \label{eq::ns_meevc_form}
		\end{equation}
		where we use the vector calculus identity $\Delta\velocity = \nabla\left(\nabla\cdot\velocity\right) - \nabla\times\nabla\times\velocity$ to derive the equality $\Delta\velocity = -\nabla\times\vorticity$. 
		
		An important aspect we wish to stress at this point is that although at the continuous level the kinematic definition $\vorticity := \nabla\times\velocity$ is valid, at the discrete level it is not always guaranteed that this identity holds. In fact, in the discretization presented here this identity is satisfied only approximately. This, as will be seen, enables the construction of a mass, energy, enstrophy and vorticity conserving discretization. 
		
	\subsection{Finite element discretization} \label{subsec::finite_element_discretization}
		In this work we set out to construct a finite element discretization for the Navier-Stokes equations as given by \eqref{eq::ns_meevc_form}. In particular we use a mixed finite element formulation, for more details on the mixed finite elements see for example \cite{brezzi1991mixed}. The first step for developing this discretization is the construction of the weak form of \eqref{eq::ns_meevc_form}:
		\begin{equation}
			\begin{dcases}
				\text{Find } \velocity\in \uspace, \totalp\in \pspace \text{ and } \vorticity \in \wspace \text{ such that:}\\
				\langle\frac{\partial\velocity}{\partial t},\ufunctionc\rangle_{\Omega} + \langle\omega\times\velocity,\ufunctionc\rangle_{\Omega} - \langle \totalp,\nabla\cdot\ufunctionc\rangle_{\Omega} = -\nu\langle\nabla\times\vorticity,\ufunctionc\rangle_{\Omega},\quad \forall \ufunctionc\in \uspace, \\
				\langle\frac{\partial\vorticity}{\partial t},\wfunctionc\rangle_{\Omega} - \frac{1}{2}\langle\vorticity,\nabla\cdot\left(\velocity\,\wfunctionc\right)\rangle_{\Omega} +  \frac{1}{2}\langle\nabla\cdot\left(\velocity\,\vorticity,\right),\wfunctionc\rangle_{\Omega} = \nu\langle\nabla\times\vorticity,\nabla\times\wfunctionc\rangle_{\Omega}, \quad\forall \wfunctionc\in \wspace, \\
				\langle\nabla\cdot\velocity,\pfunctionc\rangle_{\Omega} = 0, \quad\forall \pfunctionc\in \pspace\,,
			\end{dcases} \label{eq::ns_meevc_weak_form_continuous}
		\end{equation}
		where we have used integration by parts and the periodic boundary conditions to obtain the identities $\langle \totalp,\nabla\cdot\ufunctionc\rangle_{\Omega} = - \langle \nabla\totalp,\ufunctionc\rangle_{\Omega}$,   $\langle\vorticity,\nabla\cdot\left(\velocity\,\wfunctionc\right)\rangle_{\Omega} = -\langle\left(\velocity\cdot\nabla\right)\vorticity,\wfunctionc\rangle_{\Omega}$ and $\langle\nabla\times\vorticity,\nabla\times\wfunctionc\rangle_{\Omega} = \langle\Delta\vorticity,\wfunctionc\rangle_{\Omega}$. The space $\pspace$ corresponds to square integrable functions and the spaces $\uspace$ and $\wspace$ contain square integrable functions whose divergence and curl are also square integrable.
		
		The second step is the definition of conforming finite dimensional function spaces, where we will seek our discrete solutions for velocity $\velocityd$, pressure $\totalpd$ and vorticity $\vorticityd$:
		\begin{equation}
			\velocityd \in \udspace \subset \uspace, \quad \totalpd \in \pdspace \subset \pspace \quad \mathrm{and} \quad \vorticityd \in \wdspace \subset \wspace.
		\end{equation}
		As usual, each of these finite dimensional function spaces, $\udspace$, $\pdspace$ and $\wdspace$, has an associated finite set of basis functions, $\ubasis{i}$, $\pbasis{i}$, $\wbasis{i}$, such that
		\begin{equation}
			\udspace = \mathrm{span}\{\ubasis{1}, \dots,\ubasis{\udim}\}, \quad \pdspace = \mathrm{span}\{\pbasis{1}, \dots,\pbasis{\pdim}\}\quad\mathrm{and}\quad\wdspace = \mathrm{span}\{\wbasis{1}, \dots,\wbasis{\wdim}\},
		\end{equation}
		 where $\udim$, $\pdim$ and $\wdim$ denote the dimension of the discrete function spaces and therefore correspond to the number of degrees of freedom for each of the unknowns. As a consequence, the approximate solutions for velocity, pressure and vorticity can be expressed as a linear combination of these basis functions
		 \begin{equation}
		 	\velocityd := \sum_{i=1}^{\udim}\velocitydc{i}\,\ubasis{i}, \quad \totalpd := \sum_{i=1}^{\pdim}\totalpdc{i}\,\pbasis{i} \quad\mathrm{and}\quad \vorticityd := \sum_{i=1}^{\wdim}\vorticitydc{i}\,\wbasis{i}, \label{eq:basis_expansion}
		 \end{equation}
		 with $\velocitydc{i}$, $\totalpdc{i}$ and $\vorticitydc{i}$ the degrees of freedom of velocity, total pressure and vorticity, respectively. Since the Navier-Stokes equations form a time dependent set of equations, in general these coefficients will be time dependent, $\velocitydc{i} = \velocitydc{i}(t)$, $\totalpdc{i}=\totalpdc{i}(t)$ and $\vorticitydc{i}=\vorticitydc{i}(t)$.
		 
		The choice of the finite dimensional function spaces dictates the properties of the discretization. In order to have exact conservation of mass we must exactly satisfy the divergence free constraint at the discrete level. A sufficient condition that guarantees divergence-free velocities, e.g. \cite{arnold2010finite,Cockburn2006,Buffa2011}, is:
		\begin{equation}
			\left\{\nabla\cdot\velocityd \,|\, \velocityd\in\udspace\right\} \subseteq \pdspace\,. \label{eq:divu_subspace_q}
		\end{equation}
		In other words, the divergence operator must map $\udspace$ into $\pdspace$:
		\begin{equation}
			\udspace \stackrel{\nabla\cdot}{\longrightarrow}\pdspace\,. \label{eq:div_subcomplex}
		\end{equation}
		If we set
		\begin{equation}
			\udspace = \rtspace \quad \mathrm{and} \quad \pdspace=\qspace\,,
		\end{equation}
		 where $\rtspace$ are the Raviart-Thomas elements of degree $N$, see \cite{RaviartThomas1977,kirby2012}, and $\qspace$ are the discontinuous Lagrange elements of degree $(N-1)$, see \cite{kirby2012}, we satisfy \eqref{eq:div_subcomplex}, see for example \cite{arnold2010finite} for a proof. Additionally, this pair of finite elements satisfies the LBB stability condition, see for example \cite{RaviartThomas1977,arnold2010finite}.
		 
		What remains to define is the finite element space associated to the vorticity, $\wdspace$. Since we wish to exactly represent the diffusion term in the momentum equation $\langle\nabla\times\vorticity,\ufunction\rangle_{\Omega}$ the space $\wdspace$ must satisfy a relation similar to \eqref{eq:divu_subspace_q}
		\begin{equation}
			\left\{\nabla\times\vorticityd \,|\, \vorticityd\in\wdspace\right\} \subseteq \udspace\,. \label{eq:curlw_subspace_u}
		\end{equation}
		In other words, the curl operator must map $\wdspace$ into $\udspace$:
		\begin{equation}
			\wdspace \stackrel{\nabla\times}{\longrightarrow}\udspace\,. \label{eq:curl_subcomplex}
		\end{equation}
		The Lagrange elements of degree $N$, denoted by $\cgspace$ in \cite{kirby2012}, satisfy \eqref{eq:curl_subcomplex}, see \cite{arnold2010finite}, therefore we set
		\begin{equation}
			\wdspace = \cgspace\,.
		\end{equation}
		
		Together, the combination of these three finite element spaces forms a Hilbert subcomplex
		\begin{equation}
			0 \longrightarrow \wdspace \stackrel{\nabla\times}{\longrightarrow}\udspace \stackrel{\nabla\cdot}{\longrightarrow}\pdspace \longrightarrow 0\,,
		\end{equation}
		that mimics the 2D Hilbert complex associated to the continuous functional spaces:
		\begin{equation}
			0 \longrightarrow \wspace \stackrel{\nabla\times}{\longrightarrow}\uspace \stackrel{\nabla\cdot}{\longrightarrow}\pspace \longrightarrow 0\,.
		\end{equation}
		The Hilbert complex is an important structure that is intimately related to the de Rham complex of differential forms. The construction of a discrete subcomplex is an important requirement to obtain stable and accurate finite element discretizations, see for example \cite{arnold2010finite,Palha2014,{bossavit_japan_computational_1,bossavit_japan_computational_2,bossavit_japan_computational_3,bossavit_japan_computational_4,bossavit_japan_computational_5}} for a detailed discussion.
		
		The finite element spatial discretization used in this work to construct an approximate solution of \eqref{eq::ns_meevc_form} is then obtained by the following weak formulation
		\begin{equation}
			\begin{dcases}
				\text{Find } \velocityd\in \rtspace, \totalpd\in \qspace \text{ and } \vorticityd \in \cgspace \text{ such that:}\\
				\langle\frac{\partial\velocityd}{\partial t},\ufunction\rangle_{\Omega} + \langle\vorticityd\times\velocityd,\ufunction\rangle_{\Omega} - \langle \totalpd,\nabla\cdot\ufunction\rangle_{\Omega} = -\nu\langle\nabla\times\vorticityd,\ufunction\rangle_{\Omega},\quad \forall \ufunction\in \rtspace, \\
				\langle\frac{\partial\vorticityd}{\partial t},\wfunction\rangle_{\Omega} - \frac{1}{2}\langle\vorticityd,\nabla\cdot\left(\velocityd\,\wfunction\right)\rangle_{\Omega} + \frac{1}{2}\langle\nabla\cdot\left(\velocityd\,\vorticityd\right),\wfunction\rangle_{\Omega} = \nu\langle\nabla\times\vorticityd,\nabla\times\wfunction\rangle_{\Omega}, \quad\forall \wfunction\in \cgspace, \\
				\langle\nabla\cdot\velocityd,\pfunction\rangle_{\Omega} = 0, \quad\forall \pfunction\in \qspace\,.
			\end{dcases} \label{eq::ns_meevc_weak_form_discrete}
		\end{equation}
		
		Using the expansions for $\velocityd$, $\totalpd$ and $\vorticityd$ in \eqref{eq:basis_expansion}, \eqref{eq::ns_meevc_weak_form_discrete} can be rewritten as
		\begin{equation}
			\begin{dcases}
				\text{Find } \velocitydv\in \mathbb{R}^{\udim}, \totalpdv\in \mathbb{R}^{\pdim} \text{ and } \vorticitydv \in \mathbb{R}^{\wdim} \text{ such that:}\\
				\sum_{i=1}^{\udim}\frac{\mathrm{d}\velocitydc{i}}{\mathrm{d}t}\langle\ubasis{i},\ubasis{j}\rangle_{\Omega} + \sum_{i=1}^{\udim}\velocitydc{i}\langle\vorticityd\times\ubasis{i},\ubasis{j}\rangle_{\Omega} - \sum_{k=1}^{\pdim}\totalpdc{k}\langle \pbasis{k},\nabla\cdot\ubasis{j}\rangle_{\Omega} = -\nu\langle\nabla\times\vorticityd,\ubasis{j}\rangle_{\Omega},\quad j=1,\dots,\udim, \\
				\sum_{i=1}^{\wdim}\frac{\mathrm{d}\vorticitydc{i}}{\mathrm{d}t}\langle\wbasis{i},\wbasis{j}\rangle_{\Omega} - \sum_{i=1}^{\wdim}\frac{\vorticitydc{i}}{2}\langle\wbasis{i},\nabla\cdot\left(\velocityd\,\wbasis{j}\right)\rangle_{\Omega} + \sum_{i=1}^{\wdim}\frac{\vorticitydc{i}}{2}\langle\nabla\cdot\left(\velocityd\,\wbasis{i}\right),\wbasis{j}\rangle_{\Omega} = \nu\sum_{i=1}^{\wdim}\vorticitydc{i}\langle\nabla\times\wbasis{i},\nabla\times\wbasis{j}\rangle_{\Omega}, \quad j=1,\dots,\wdim, \\
				\sum_{i=1}^{\udim}\velocitydc{i}\langle\nabla\cdot\ubasis{i},\pbasis{j}\rangle_{\Omega} = 0, \quad j = 1,\dots,\pdim\,,
			\end{dcases} \label{eq::ns_meevc_weak_form_discrete_expansion}
		\end{equation}
		with $\velocitydv := [\velocitydc{1},\dots,\velocitydc{\udim}]^{\top}$, $\totalpdv := [\totalpdc{1},\dots,\totalpdc{\pdim}]^{\top}$ and $\vorticitydv := [\vorticitydc{1},\dots,\vorticitydc{\wdim}]^{\top}$. Using matrix notation, \eqref{eq::ns_meevc_weak_form_discrete_expansion} can be expressed more compactly as
		\begin{equation}
			\begin{dcases}
				\text{Find } \velocitydv\in \mathbb{R}^{\udim}, \totalpdv\in \mathbb{R}^{\pdim} \text{ and } \vorticitydv \in \mathbb{R}^{\wdim} \text{ such that:}\\
				\matrixoperator{M} \frac{\mathrm{d}\velocitydv}{\mathrm{d}t} + \matrixoperator{R}\,\velocitydv - \matrixoperator{P}\,\totalpdv = -\nu\,\boldsymbol{l}, \\
				\matrixoperator{N}\frac{\mathrm{d}\vorticitydv}{\mathrm{d}t}  - \frac{1}{2}\matrixoperator{W}\,\vorticitydv + \frac{1}{2}\matrixoperator{W}^{\top}\vorticitydv = \nu\,\matrixoperator{L}\,\vorticitydv, \\
				\matrixoperator{D}\,\velocitydv = 0,
			\end{dcases} \label{eq::ns_meevc_weak_form_discrete_matrix_notation}
		\end{equation}
		The coefficients of the matrices $\matrixoperator{M}$, $\matrixoperator{R}$ and $\matrixoperator{P}$, and the column vector $\boldsymbol{l}$ are given by
		\begin{equation}
			\matrixoperatorc{M}_{ij} := \langle\ubasis{j},\ubasis{i}\rangle_{\Omega}, \quad \matrixoperatorc{R}_{ij} := \langle\vorticityd\times\ubasis{j},\ubasis{i}\rangle_{\Omega}, \quad \matrixoperatorc{P}_{ij} := \langle \pbasis{j},\nabla\cdot\ubasis{i}\rangle_{\Omega}\quad\mathrm{and}\quad l_{i} := \langle\nabla\times\vorticityd,\ubasis{i}\rangle_{\Omega}. \label{eq:matrix_coefficients_1}
		\end{equation}
		Similarly, the coefficients of the matrices $\matrixoperator{N}$, $\matrixoperator{W}$, $\matrixoperator{L}$ and $\matrixoperator{D}$ are given by
		\begin{equation}
			\matrixoperatorc{N}_{ij} := \langle\wbasis{j},\wbasis{i}\rangle_{\Omega}, \quad \matrixoperatorc{W}_{ij} := \langle\wbasis{j},\nabla\cdot\left(\velocityd\,\wbasis{i}\right)\rangle_{\Omega}, \quad \matrixoperator{L}_{ij} := \langle\nabla\times\wbasis{j},\nabla\times\wbasis{i}\rangle_{\Omega} \quad \mathrm{and}\quad  \matrixoperator{D}_{ij} :=\langle\nabla\cdot\ubasis{j},\pbasis{i}\rangle_{\Omega}. \label{eq:matrix_coefficients_2}
		\end{equation}

\section{Temporal discretization} \label{sec::temporal_discretization}
	In this section we present the time discretization. The choice of a time discretization is essential in preserving invariants. Not all time integrators satisfy conservation of energy, even though the spatial discretization imposes conservation of kinetic energy as a function of time, see \cite{Sanderse2013}. In this work we choose to use a Gauss method. This time integrator is a type of collocation method based on Gauss quadrature. It is known to  be an implicit Runge-Kutta method and to have optimal convergence order $2s$ for $s$ stages. Two of its most attractive properties are that (i) it conserves energy when applied to the discretization of the Navier-Stokes equations, see \cite{Sanderse2013}, and (ii) it is time-reversible, see \cite{Hairer2006}. For a more detailed discussion of its properties and construction see \cite{Hairer2006}. Although any Gauss integrator could be used, we choose the lowest order, $s=1$, also known as the \emph{midpoint rule}, because it enables the construction of an explicit staggered integrator in time.

	When applied to the solution of a 1D ordinary differential equation of the form
	\begin{equation}
		\begin{dcases}
			\frac{\mathrm{d}f}{\mathrm{d}t} = g(f(t),t), \\
			f(0) = f_{0},
		\end{dcases}
	\end{equation}
	the one stage Gauss integrator results in the following implicit time stepping scheme
	\begin{equation}
		\frac{f^{k} - f^{k-1}}{\Delta t} = g\left(\frac{f^{k}+f^{k-1}}{2},t+\frac{\Delta t}{2}\right), \quad k=1,\dots,M, \label{eq:gauss_integrator_1D}
	\end{equation}
	where $f^{0} = f_{0}$, $\Delta t$ is the time step and $M$ is the number of time steps. The direct application of \eqref{eq:gauss_integrator_1D} to the discrete weak form \eqref{eq::ns_meevc_weak_form_discrete} results in
	\begin{equation}
			\begin{dcases}
				\text{Find } \velocitydnext\in \rtspace, \totalpdnext\in \qspace \text{ and } \vorticitydnext \in \cgspace \text{ such that:}\\
				\langle\frac{\velocitydnext - \velocitydprevious}{\Delta t},\ufunction\rangle_{\Omega} + \langle\vorticitydmida\times\frac{\velocitydnext+\velocitydprevious}{2},\ufunction\rangle_{\Omega} - \langle \totalpdmida,\nabla\cdot\ufunction\rangle_{\Omega} = -\nu\langle\nabla\times\vorticitydmida,\ufunction\rangle_{\Omega},\quad \forall \ufunction\in \rtspace, \\
				\langle\frac{\vorticitydnext-\vorticitydprevious}{\Delta t},\wfunction\rangle_{\Omega} - \frac{1}{2}\langle\frac{\vorticitydnext+\vorticitydprevious}{2},\nabla\cdot\left(\velocitydmida\,\wfunction\right)\rangle_{\Omega} + \frac{1}{2}\langle\nabla\cdot\left(\velocitydmida\,\frac{\vorticitydnext+\vorticitydprevious}{2}\right),\wfunction\rangle_{\Omega}  = \nu\langle\nabla\times\frac{\vorticitydnext+\vorticitydprevious}{2},\nabla\times\wfunction\rangle_{\Omega}, \quad\forall \wfunction\in \cgspace, \\
				\langle\nabla\cdot\velocitydnext,\pfunction\rangle_{\Omega} = 0, \quad\forall \pfunction\in \qspace\,,
			\end{dcases} \label{eq::ns_meevc_weak_form_discrete_naive_gauss}
		\end{equation}
	where, for compactness of notation, we have set
	\begin{equation}
		\velocitydmida := \frac{\velocitydnext + \velocitydprevious}{2} \quad \mathrm{and}\quad \vorticitydmida := \frac{\vorticitydnext + \vorticitydprevious}{2}. \label{eq:mid_steps_compact_notation}
	\end{equation}
	
\begin{figure}[!ht]
		\centering
		\includegraphics{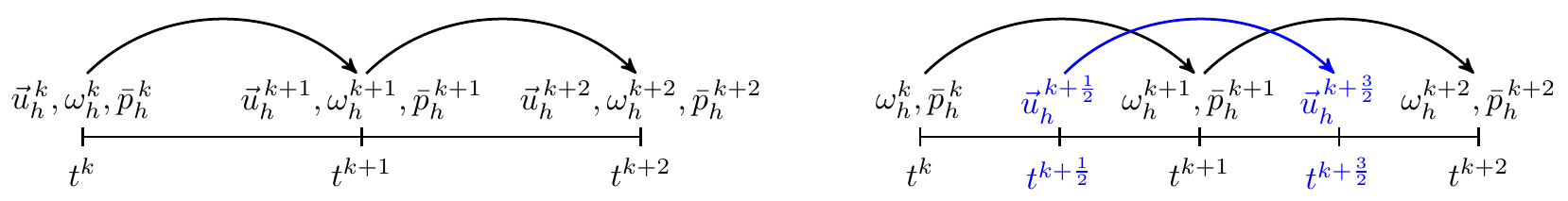}
		\caption{Diagram of the time stepping. Left: all unknowns at the same time instant, as in \eqref{eq::ns_meevc_weak_form_discrete_naive_gauss}. Right: staggered in time unknowns, as in \eqref{eq::ns_meevc_weak_form_discrete_staggered_gauss}.}
		\label{fig:time_stepping}
	\end{figure}

	The time stepping scheme in \eqref{eq::ns_meevc_weak_form_discrete_naive_gauss} consists of a coupled system of nonlinear equations. Therefore its solution will necessarily require an iterative procedure, which is computationally expensive. To overcome this drawback, instead of defining all the unknown physical quantities $\velocityd$, $\vorticityd$ and  $\totalpd$, at the same time instants $t^{k}$ we choose to stagger them in time. In this way it is possible to transform \eqref{eq::ns_meevc_weak_form_discrete_naive_gauss} into two systems of quasi-linear equations. The unknown vorticity and total pressure are defined at the integer time instants $\vorticitydprevious$, $\totalpdprevious$ and the unknown velocity is defined at the intermediate time instants $\velocitydmidb$, see \figref{fig:time_stepping}. Taking into account this staggered approach, \eqref{eq::ns_meevc_weak_form_discrete_naive_gauss} can be rewritten as
	\begin{equation}
			\begin{dcases}
				\text{Find } \velocitydmidnext\in \rtspace, \totalpdnext\in \qspace \text{ and } \vorticitydnext \in \cgspace \text{ such that:}\\
				\langle\frac{\velocitydmidnext - \velocitydmidb}{\Delta t},\ufunction\rangle_{\Omega} + \langle\vorticitydnext\times\frac{\velocitydmidnext+\velocitydmidb}{2},\ufunction\rangle_{\Omega} - \langle \totalpdnext,\nabla\cdot\ufunction\rangle_{\Omega} = -\nu\langle\nabla\times\vorticitydnext,\ufunction\rangle_{\Omega},\quad \forall \ufunction\in \rtspace, \\
				\langle\frac{\vorticitydnext-\vorticitydprevious}{\Delta t},\wfunction\rangle_{\Omega}  - \frac{1}{2}\langle\frac{\vorticitydnext+\vorticitydprevious}{2},\nabla\cdot\left(\velocitydmidb\,\wfunction\right)\rangle_{\Omega} + \frac{1}{2}\langle\nabla\cdot\left(\velocitydmidb\,\frac{\vorticitydnext+\vorticitydprevious}{2}\right),\wfunction\rangle_{\Omega} = \nu\langle\nabla\times\frac{\vorticitydnext+\vorticitydprevious}{2},\nabla\times\wfunction\rangle_{\Omega}, \quad\forall \wfunction\in \cgspace, \\
				\langle\nabla\cdot\velocitydmidnext,\pfunction\rangle_{\Omega} = 0, \quad\forall \pfunction\in \qspace\,,
			\end{dcases} \label{eq::ns_meevc_weak_form_discrete_staggered_gauss}
		\end{equation}
	where $\velocitydmidb$ and $\vorticitydprevious$ are known at the start of each time step. 
	
	Using \eqref{eq::ns_meevc_weak_form_discrete_matrix_notation}, it is possible to rewrite \eqref{eq::ns_meevc_weak_form_discrete_staggered_gauss} in a compact matrix notation
	\begin{equation}
			\begin{dcases}
				\text{Find } \velocitydvmidnext\in \mathbb{R}^{\udim}, \totalpdvnext\in \mathbb{R}^{\pdim} \text{ and } \vorticitydvnext \in \mathbb{R}^{\wdim} \text{ such that:}\\
				\matrixoperator{M} \frac{\velocitydvmidnext - \velocitydvmidb}{\Delta t} + \matrixoperator{R}^{k+1}\,\frac{\velocitydvmidnext + \velocitydvmidb}{2} - \matrixoperator{P}\,\totalpdvnext = -\nu\,\boldsymbol{l}^{k+1}, \\
				\matrixoperator{N}\frac{\vorticitydvnext - \vorticitydvprevious}{\Delta t}  - \frac{1}{2}\matrixoperator{W}^{k+\frac{1}{2}}\,\frac{\vorticitydvnext + \vorticitydvprevious}{2} + \frac{1}{2}\left(\matrixoperator{W}^{k+\frac{1}{2}}\right)^{\top}\frac{\vorticitydvnext + \vorticitydvprevious}{2} = \nu\,\matrixoperator{L}\,\frac{\vorticitydvnext + \vorticitydvprevious}{2}, \\
				\matrixoperator{D}\,\velocitydvmidnext = 0,
			\end{dcases} \label{eq::ns_meevc_weak_form_discrete_staggered_gauss_matrix_notation}
	\end{equation}
	where all matrix operators are as in \eqref{eq:matrix_coefficients_1} and \eqref{eq:matrix_coefficients_2}, with the exception of $\matrixoperator{R}^{k+1}$, $\matrixoperator{W}^{k+\frac{1}{2}}$ and $\boldsymbol{l}^{k+1}$, the coefficients of which are
	\begin{equation}
		\matrixoperatorc{R}^{k+1}_{ij} := \langle\vorticitydnext\times\ubasis{j},\ubasis{i}\rangle_{\Omega}, \quad l^{k+1}_{i} := \langle\nabla\times\vorticitydnext,\ubasis{i}\rangle_{\Omega}\quad\mathrm{and}\quad \matrixoperatorc{W}^{k+\frac{1}{2}}_{ij} := \langle\wbasis{j},\nabla\cdot\left(\velocitydmidb\,\wbasis{i}\right)\rangle_{\Omega}. \label{eq:matrix_coefficients_staggered}
	\end{equation}
	
	To start the iteration procedure $\velocitydmidstart$ and $\vorticitydstart$ are required. Since only $\velocitydstart$ and $\vorticitydstart$ are known, the first time step needs to be implicit, according to \eqref{eq::ns_meevc_weak_form_discrete_naive_gauss}. Once $\velocitydstartnext$ is known, \eqref{eq:mid_steps_compact_notation} can be used to retrieve $\velocitydmidstart$. The remaining time steps can then be computed explicitly with \eqref{eq::ns_meevc_weak_form_discrete_staggered_gauss}.
	
\section{Conservation properties and time reversibility} \label{sec::conservation_properties_and_time_reversibility}
	It is known that in the absence of external forces and viscosity the incompressible Navier-Stokes equations contain important symmetries or invariants, e.g. \cite{Arnold1966,Arnold1992,Majda2001,Foias2001}. Conservation of mass, energy, enstrophy and total vorticity are such invariants.
	
	A discretization for the incompressible Navier-Stokes equations has been proposed in \eqref{eq::ns_meevc_weak_form_discrete_staggered_gauss}. In the following sections we prove its conservation properties, with respect to mass, energy, enstrophy and vorticity. Additionally, we also prove time reversibility.
	
	\subsection{Mass and vorticity conservation} \label{sec:mass_vorticity_conservation}
		Mass conservation is given by the divergence of the velocity field. Since with this discretization the discrete velocity field is exactly divergence free mass is conserved.
		
		Regarding vorticity conservation, the time evolution of vorticity is governed by the vorticity transport equation, \eqref{eq::ns_meevc_form}:
		\begin{equation}
			\frac{\partial \vorticity}{\partial t} +\frac{1}{2}\left(\velocity\cdot\nabla\right)\vorticity + \frac{1}{2}\nabla\cdot\left(\velocity\vorticity\right) = \nu \Delta\vorticity .
		\end{equation}
		This equation, by \eqref{eq::ns_meevc_weak_form_discrete_staggered_gauss}, is discretized as
		\begin{equation}
			\begin{dcases}
				\text{Find } \gentranspdnext \in \cgspace \text{ such that:}\\
				\langle\frac{\gentranspdnext-\gentranspdprevious}{\Delta t},\wfunction\rangle_{\Omega}  - \frac{1}{2}\langle\frac{\gentranspdnext+\gentranspdprevious}{2},\nabla\cdot\left(\velocitydmidb\,\wfunction\right)\rangle_{\Omega} + \frac{1}{2}\langle\nabla\cdot\left(\velocitydmidb\,\frac{\gentranspdnext+\gentranspdprevious}{2}\right),\wfunction\rangle_{\Omega} = \nu\langle\nabla\times\frac{\gentranspdnext+\gentranspdprevious}{2},\nabla\times\wfunction\rangle_{\Omega}, \quad\forall \wfunction\in \cgspace,
			\end{dcases} \label{eq::generic_transport_equation}
		\end{equation}
		Conservation of vorticity $\gentranspd$ corresponds to
		\begin{equation}
			\int_{\Omega}\gentranspdnext := \langle\gentranspdnext,1\rangle =  \langle\gentranspdprevious,1\rangle := \int_{\Omega}\gentranspdprevious.
		\end{equation}
		Therefore, the proof of conservation of vorticity is obtained from evaluating \eqref{eq::generic_transport_equation} for the special case $\wfunction = 1$:
		\begin{equation}
			\langle\frac{\gentranspdnext-\gentranspdprevious}{\Delta t},1\rangle_{\Omega}  - \frac{1}{2}\langle\frac{\gentranspdnext+\gentranspdprevious}{2},\nabla\cdot\velocitydmidb\rangle_{\Omega} + \frac{1}{2}\langle\nabla\cdot\left(\velocitydmidb\,\frac{\gentranspdnext+\gentranspdprevious}{2}\right),1\rangle_{\Omega} = \nu\langle\nabla\times\frac{\gentranspdnext+\gentranspdprevious}{2},\nabla\times1\rangle_{\Omega} . \label{eq:proof_conservation_mass_vorticity_1}
		\end{equation}
		Since the function spaces have been chosen such that $\nabla\cdot\velocityd = 0$ is satisfied exactly, see \secref{subsec::finite_element_discretization}, equation \eqref{eq:proof_conservation_mass_vorticity_1} becomes
		\begin{equation}
			\langle\gentranspdnext,1\rangle_{\Omega}  + \frac{\Delta t}{2}\langle\nabla\cdot\left(\velocitydmidb\,\frac{\gentranspdnext+\gentranspdprevious}{2}\right),1\rangle_{\Omega} = \langle\gentranspdprevious,1\rangle_{\Omega} . \label{eq:proof_conservation_mass_vorticity_2}
		\end{equation}
		The second term on the left hand side can be rewritten as
		\begin{equation}
			\langle\nabla\cdot\left(\velocitydmidb\,\frac{\gentranspdnext+\gentranspdprevious}{2}\right),1\rangle_{\Omega} := \int_{\Omega} \nabla\cdot\left(\velocitydmidb\,\frac{\gentranspdnext+\gentranspdprevious}{2}\right) = \int_{\partial\Omega} \left(\velocitydmidb\,\frac{\gentranspdnext+\gentranspdprevious}{2}\right)\cdot\vec{n}, \label{eq:proof_conservation_mass_vorticity_3}
		\end{equation}
		with $\vec{n}$ the outward unit normal at the boundary of $\Omega$. At the continuous level, the boundary integral in  \eqref{eq:proof_conservation_mass_vorticity_3} is trivially equal to zero for periodic boundary conditions, which is the case we will consider here. Note that since the domain is periodic, the nodes of the mesh at opposite sides of domain must coincide, otherwise the periodicity is lost. Since $\gentranspd \in \cgspace$, it is continuous across elements. In a similar way, the normal component of $\velocityd$ is continuous across elements because $\velocityd \in \rtspace$, see for example \cite{RaviartThomas1977,arnold2010finite,kirby2012}. Therefore, in the case of periodic boundary conditions the boundary integral in \eqref{eq:proof_conservation_mass_vorticity_3} will still be exactly equal to zero and \eqref{eq:proof_conservation_mass_vorticity_2} becomes
		\begin{equation}
			\langle\gentranspdnext,1\rangle_{\Omega}  = \langle\gentranspdprevious,1\rangle_{\Omega},\label{eq:proof_conservation_mass_vorticity_4}
		\end{equation}
		which proves that vorticity $\vorticitytdprevious$ is conserved because
		\begin{equation}
			\vorticitytdnext := \langle\vorticityd^{k+1},1\rangle_{\Omega}  = \langle\vorticityd^{k},1\rangle_{\Omega} =:\vorticitytdprevious
		\end{equation}
		 is conserved.
		
	\subsection{Kinetic energy conservation} \label{sec::kinetic_energy_conservation}
		Kinetic energy conservation is one of the two \emph{secondary conservation} properties of the numerical solver proposed in this work. Here we prove that total kinetic energy $\energy$ is conserved by \eqref{eq::ns_meevc_weak_form_discrete_staggered_gauss}.
		
		At the continuous level, kinetic energy is defined as proportional to the $L^{2}(\Omega)$ norm of the velocity velocity
		\begin{equation}
			\energy := \frac{1}{2}\|\velocity\|_{L^{2}(\Omega)} := \frac{1}{2}\langle\velocity,\velocity\rangle,
		\end{equation}
		is conserved in the inviscid limit $\nu = 0$, see for example \cite{Majda2001,Foias2001}.
		This definition can be directly extended to the discrete level as
		\begin{equation}
			\energyd := \frac{1}{2}\|\velocityd\|_{L^{2}(\Omega)} := \frac{1}{2}\langle\velocityd,\velocityd\rangle.
		\end{equation}
		Therefore, kinetic energy conservation at the discrete level corresponds to
		\begin{equation}
			\energydnext := \frac{1}{2}\langle\velocitydmidnext,\velocitydmidnext\rangle =  \frac{1}{2}\langle\velocitydmidb,\velocitydmidb\rangle := \energydprevious. \label{eq:energy_conservation_definition}
		\end{equation}
		
		To prove this identity we take the first equation in \eqref{eq::ns_meevc_weak_form_discrete_staggered_gauss} (momentum equation) in the inviscid limit $\nu = 0$ and choose $\ufunction = \frac{1}{2}\velocitydmidnext + \frac{1}{2}\velocitydmidb$, resulting in
		\begin{equation}
			 \frac{1}{2}\langle\frac{\velocitydmidnext - \velocitydmidb}{\Delta t},\velocitydmidnext + \velocitydmidb\rangle_{\Omega} +  \frac{1}{2}\langle\vorticitydnext\times\frac{\velocitydmidnext+\velocitydmidb}{2},\velocitydmidnext + \velocitydmidb\rangle_{\Omega} -  \frac{1}{2}\langle \totalpdnext,\nabla\cdot\left(\velocitydmidnext + \velocitydmidb\right)\rangle_{\Omega} =0. \label{eq:energy_conservation_proof_1}
		\end{equation}
		The term involving the total pressure $\totalpd$ is identically zero because the velocity field is divergence free at every time step, due to the particular choice of function spaces, as discussed in \secref{subsec::finite_element_discretization}. After rearranging, equation \eqref{eq:energy_conservation_proof_1} then becomes
		 \begin{equation}
			 \frac{1}{2}\langle\velocitydmidnext,\velocitydmidnext\rangle_{\Omega}  + \frac{\Delta t}{2}\,\langle\vorticitydnext\times\frac{\velocitydmidnext+\velocitydmidb}{2},\velocitydmidnext + \velocitydmidb\rangle_{\Omega} =  \frac{1}{2}\langle\velocitydmidb,\velocitydmidb\rangle_{\Omega}. \label{eq:energy_conservation_proof_1b}
		\end{equation}
		Because of Lemma 1.3 in \cite{TemamNS}, \eqref{eq:energy_conservation_proof_1b} becomes
		\[
			 \frac{1}{2}\langle\velocitydmidnext,\velocitydmidnext\rangle_{\Omega}  =  \frac{1}{2}\langle\velocitydmidb,\velocitydmidb\rangle_{\Omega},
		\]
		which proves \eqref{eq:energy_conservation_definition}.
		
	\subsection{Enstrophy conservation} \label{sec::enstrophy_conservation}
		The second \emph{secondary conservation} property of the numerical solver presented here is enstrophy conservation. We then proceed to prove conservation of enstrophy $\enstrophy$ at the discrete level.
		
		Similarly to kinetic energy, enstrophy is defined as proportional to the $L^{2}(\Omega)$ norm of vorticity
		\begin{equation}
			\enstrophy := \frac{1}{2} \|\vorticity\|_{L^{2}(\Omega)} := \langle\vorticity,\vorticity\rangle_{\Omega},
		\end{equation}
		and is also a conserved quantity in the inviscid limit $\nu = 0$, see for example \cite{Majda2001,Foias2001}. This definition can be straightforwardly extended to the discrete level as
		\begin{equation}
			\enstrophyd := \frac{1}{2} \|\vorticityd\|_{L^{2}(\Omega)} := \langle\vorticityd,\vorticityd\rangle_{\Omega},
		\end{equation}
		This directly implies that conservation of enstrophy at the discrete level requires that
		\begin{equation}
			\enstrophydnext := \frac{1}{2}\langle\vorticitydnext,\vorticitydnext\rangle =  \frac{1}{2}\langle\vorticitydprevious,\vorticitydprevious\rangle := \enstrophydprevious.
		\end{equation}
		
		The proof of enstrophy conservation follows a similar approach to the one used for the proof of energy conservation. In this case we start with the second equation in \eqref{eq::ns_meevc_weak_form_discrete_staggered_gauss} (vorticity transport equation) in the inviscid limit $\nu = 0$ and choose $\wfunction = \frac{1}{2}\vorticitydnext + \frac{1}{2}\vorticitydprevious$, obtaining
		\begin{equation}
			\frac{1}{2}\langle\frac{\vorticitydnext-\vorticitydprevious}{\Delta t},\vorticitydnext+\vorticitydprevious\rangle_{\Omega}  - \frac{1}{4}\langle\frac{\vorticitydnext+\vorticitydprevious}{2},\nabla\cdot\left(\velocitydmidb\,\left(\vorticitydnext+\vorticitydprevious\right)\right)\rangle_{\Omega} + \frac{1}{4}\langle\nabla\cdot\left(\velocitydmidb\,\frac{\vorticitydnext+\vorticitydprevious}{2}\right),\vorticitydnext+\vorticitydprevious\rangle_{\Omega} = 0. \label{eq:enstrophy_conservation_proof_1}
		\end{equation}
		The last two terms on the left hand side of \eqref{eq:enstrophy_conservation_proof_1} are equal and have opposite signs, therefore cancelling each other. Conservation of enstrophy then follows directly.
		
	\subsection{Time reversibility}
		It is well known that the Navier-Stokes equations in the inviscid limit $\nu = 0$ are time reversible, e.g. \cite{Majda2001,Duponcheel2008}. This important property has served in the past as a benchmark for numerical flow solvers \cite{Duponcheel2008} and for the construction of improved LES subgrid models \cite{CARATI2001}, for example. It is intimately related to the numerical dissipation and therefore it is desirable to satisfy it.
		
		To prove time reversibility of the proposed scheme we follow the same approach presented by Hairer et al. in \cite{Hairer2006}. Our time integrator $\Phi_{\Delta t}$ is a \emph{one-step} method since it uses only the information of the previous time step:
		\[
			\Phi_{\Delta t} (\velocitydmidb,\vorticitydprevious) =  (\velocitydmidnext,\vorticitydnext).
		\]
		A numerical one-step method $\Phi_{\Delta t}$ is \emph{time-reversible}, if it satisfies, see Hairer et al. \cite{Hairer2006},
		\[
			\Phi_{\Delta t}\circ\Phi_{-\Delta t} = id \qquad \text{or equivalently} \qquad \Phi_{\Delta t} = \Phi^{-1}_{-\Delta t}\,.
		\]
		To prove time reversibility of our method we will show that $\Phi_{\Delta t}\circ\Phi_{-\Delta t} = id$.
		
\begin{figure}[!ht]
			\centering
			\includegraphics{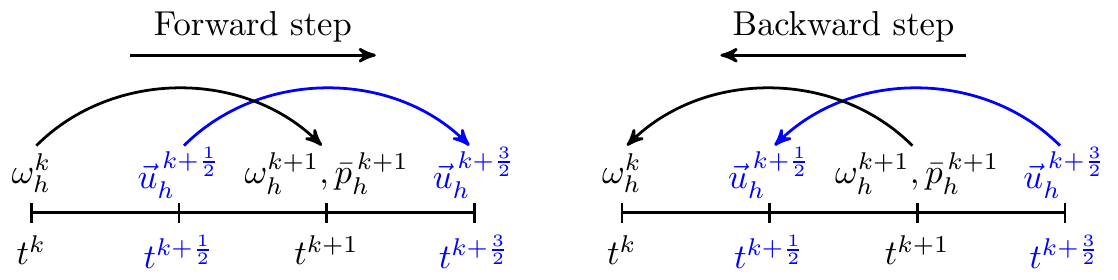}
			\caption{Diagram of time reversibility. Left: forward step. Right: backward step.}
			\label{fig:time_reversibility}
		\end{figure}

		Start with known velocity $\velocitydreverseprevious$ at time instant $t^{k+\frac{1}{2}}$ and known vorticity $\vorticitydreverseprevious$ at time instant $t^{k}$. Advance both quantities one time step to obtain the velocity $\velocitydreversenext$ at time instant $t^{k+\frac{3}{2}}$ and the vorticity $\vorticitydreversenext$ at the time instant $t^{k+1}$. Reverse time and compute one time step. If the numerical method is time reversible the initial velocity and vorticity fields must be retrieved, as represented in \figref{fig:time_reversibility}. First we prove the reversibility of the vorticity time step and then the reversibility of the velocity time step.

		\subsubsection{Reversibility of vorticity time step}
			In the inviscid limit $\nu = 0$, the forward vorticity time step, \figref{fig:time_reversibility} left, which advances the vorticity field from the time instant $t^{k}$ to the time instant $t^{k+1}$ is given by the second equation in \eqref{eq::ns_meevc_weak_form_discrete_staggered_gauss} (vorticity transport equation)
			\begin{equation}
				\begin{dcases}
					\text{Find } \vorticitydnext \in \cgspace \text{ such that:}\\
					\langle\frac{\vorticitydnext-\vorticitydprevious}{\Delta t},\wfunction\rangle_{\Omega}  - \frac{1}{2}\langle\frac{\vorticitydnext+\vorticitydprevious}{2},\nabla\cdot\left(\velocitydmidb\,\wfunction\right)\rangle_{\Omega} + \frac{1}{2}\langle\nabla\cdot\left(\velocitydmidb\,\frac{\vorticitydnext+\vorticitydprevious}{2}\right),\wfunction\rangle_{\Omega} = 0,\quad\forall\wfunction\in\cgspace. 
				\end{dcases}\label{eq:time_reversibility_vorticity_proof_1}
			\end{equation}
			Using \eqref{eq::ns_meevc_weak_form_discrete_staggered_gauss_matrix_notation}, this expression  can be rewritten in matrix notation as an algebraic system of equations:
			\begin{equation}
				\matrixoperator{N}\frac{\vorticitydvnext - \vorticitydvprevious}{\Delta t}  - \frac{1}{2}\matrixoperator{W}^{k+\frac{1}{2}}\,\frac{\vorticitydvnext + \vorticitydvprevious}{2} + \frac{1}{2}\left(\matrixoperator{W}^{k+\frac{1}{2}}\right)^{\top}\frac{\vorticitydvnext + \vorticitydvprevious}{2} = 0.\label{eq:time_reversibility_vorticity_proof_1b}
			\end{equation}
			Rearranging this expression it is possible to solve for $\vorticitydvnext$
			\begin{equation}
				\vorticitydvnext = \left(\matrixoperator{N}  - \frac{\Delta t}{4}\matrixoperator{W}^{k+\frac{1}{2}} + \frac{\Delta t}{4}\left(\matrixoperator{W}^{k+\frac{1}{2}}\right)^{\top}\right)^{-1}\left(\matrixoperator{N}  + \frac{\Delta t}{4}\matrixoperator{W}^{k+\frac{1}{2}} - \frac{\Delta t}{4}\left(\matrixoperator{W}^{k+\frac{1}{2}}\right)^{\top}\right)\vorticitydvprevious . \label{eq:time_reversibility_vorticity_proof_2}
			\end{equation}
			If we now introduce the compact notation
			\begin{equation}
				\matrixoperator{A}_{-}:= \matrixoperator{N}  - \frac{\Delta t}{4}\matrixoperator{W}^{k+\frac{1}{2}} + \frac{\Delta t}{4}\left(\matrixoperator{W}^{k+\frac{1}{2}}\right)^{\top} \quad\mathrm{and}\quad\matrixoperator{A}_{+}:=\matrixoperator{N}  + \frac{\Delta t}{4}\matrixoperator{W}^{k+\frac{1}{2}} - \frac{\Delta t}{4}\left(\matrixoperator{W}^{k+\frac{1}{2}}\right)^{\top},
			\end{equation}
			equation \eqref{eq:time_reversibility_vorticity_proof_2} becomes
			\begin{equation}
				\vorticitydvnext = \matrixoperator{A}_{-}^{-1}\matrixoperator{A}_{+}\,\vorticitydvprevious. \label{eq:vorticity_reversibility_forward_compact}
			\end{equation}
			
			The backward time step for vorticity, \figref{fig:time_reversibility} right, is obtained by substituting $\Delta t$ by $-\Delta t$. Substituting this transformation in \eqref{eq:time_reversibility_vorticity_proof_1} and using \eqref{eq::ns_meevc_weak_form_discrete_staggered_gauss_matrix_notation} leads to an expression identical to \eqref{eq:time_reversibility_vorticity_proof_1b}
			\begin{equation}
				\matrixoperator{N}\frac{\vorticitydvprevious - \vorticitydvnext}{-\Delta t}  - \frac{1}{2}\matrixoperator{W}^{k+\frac{1}{2}}\,\frac{\vorticitydvprevious + \vorticitydvnext}{2} + \frac{1}{2}\left(\matrixoperator{W}^{k+\frac{1}{2}}\right)^{\top}\frac{\vorticitydvprevious + \vorticitydvnext}{2} = 0.\label{eq:time_reversibility_vorticity_proof_3}
			\end{equation}
			This expression can be rearranged to yield a result similar to \eqref{eq:vorticity_reversibility_forward_compact}
			\begin{equation}
				\vorticitydvprevious = \matrixoperator{A}_{+}^{-1}\matrixoperator{A}_{-}\,\vorticitydvnext. \label{eq:vorticity_reversibility_forward_compact_2}
			\end{equation}
			Combining the forward time step \eqref{eq:vorticity_reversibility_forward_compact} with the backward time step	\eqref{eq:vorticity_reversibility_forward_compact_2} results in
			\begin{equation}
				\vorticitydvprevious = \matrixoperator{A}_{+}^{-1}\matrixoperator{A}_{-}\left(\matrixoperator{A}_{-}^{-1}\matrixoperator{A}_{+}\vorticitydvprevious\right) = \vorticitydvprevious.
			\end{equation}
			Thus showing the reversibility of the vorticity time step.
			
		\subsubsection{Reversibility of velocity time step}
			For the reversibility of the velocity time step consider first the forward time step, \figref{fig:time_reversibility} left, in the inviscid limit $\nu = 0$. The first equation in \eqref{eq::ns_meevc_weak_form_discrete_staggered_gauss} (momentum equation) computes the evolution of the velocity field from time instant $t^{k+\frac{1}{2}}$ to the time instant $t^{k+\frac{3}{2}}$
				\begin{equation}
				\begin{dcases}
					\text{Find } \velocitydmidnext\in \rtspace, \totalpdnext\in \qspace \text{ such that:}\\
				  		\langle\frac{\velocitydmidnext - \velocitydmidb}{\Delta t},\ufunction\rangle_{\Omega} + \langle\vorticitydnext\times\frac{\velocitydmidnext+\velocitydmidb}{2},\ufunction\rangle_{\Omega} - \langle \totalpdnext,\nabla\cdot\ufunction\rangle_{\Omega} = 0,\quad \forall \ufunction\in \rtspace, \\
				\langle\nabla\cdot\velocitydmidnext,\pfunction\rangle_{\Omega} = 0, \quad\forall \pfunction\in \qspace\,.
				\end{dcases} \label{eq::velocity_reversibility_proof_1}
			\end{equation}				
			Using \eqref{eq::ns_meevc_weak_form_discrete_staggered_gauss_matrix_notation} we can write this expression as an algebraic system of equations
			\begin{equation}
				\begin{dcases}
					\matrixoperator{M} \frac{\velocitydvmidnext - \velocitydvmidb}{\Delta t} + \matrixoperator{R}^{k+1}\,\frac{\velocitydvmidnext + \velocitydvmidb}{2} - \matrixoperator{P}\,\totalpdvnext =0, \\
					\matrixoperator{D}\,\velocitydvmidnext = 0,
				\end{dcases} \label{eq::velocity_reversibility_proof_2}
			\end{equation}
			of which $\velocitydvmidnext$ and $\totalpdvnext$ are the solution. Once $\totalpdvnext$ is known, it is possible to write an explicit expression for $\velocitydvmidnext$ as a function of $\velocitydvmidb$ and $\totalpdvnext$
			\begin{equation}
				\velocitydvmidnext = \left(\matrixoperator{M} +\frac{\Delta t}{2} \matrixoperator{R}^{k+1}\right)^{-1}\left(\matrixoperator{M} -\frac{\Delta t}{2} \matrixoperator{R}^{k+1}\right) \velocitydvmidb+ \Delta t \,\left(\matrixoperator{M} +\frac{\Delta t}{2} \matrixoperator{R}^{k+1}\right)^{-1}\matrixoperator{P}\,\totalpdvnext.  \label{eq:velocity_reversibility_proof_3}
			\end{equation}
			Introducing the compact notation
			\begin{equation}
				\matrixoperator{B}_{+} := \matrixoperator{M} +\frac{\Delta t}{2} \matrixoperator{R}^{k+1} \quad\mathrm{and}\quad \matrixoperator{B}_{-} := \matrixoperator{M} -\frac{\Delta t}{2} \matrixoperator{R}^{k+1},
			\end{equation}
			equation \eqref{eq:velocity_reversibility_proof_3} becomes
			\begin{equation}
				\velocitydvmidnext = \matrixoperator{B}_{+}^{-1}\matrixoperator{B}_{-}\,\velocitydvmidb + \Delta t \,\matrixoperator{B}_{+}^{-1}\matrixoperator{P}\,\totalpdvnext.  \label{eq:velocity_reversibility_proof_4}
			\end{equation}
			
			To compute the backward time step for velocity, \figref{fig:time_reversibility} right, we proceed in the same manner as for vorticity: reverse the time step by replacing $\Delta t$ by $-\Delta t$. We can now apply this transformation to \eqref{eq::velocity_reversibility_proof_1} and use \eqref{eq::ns_meevc_weak_form_discrete_staggered_gauss_matrix_notation}  to obtain an expression identical to \eqref{eq::velocity_reversibility_proof_2}
			\begin{equation}
				\begin{dcases}
					\matrixoperator{M} \frac{\velocitydvmidnext - \velocitydvmidb}{-\Delta t} + \matrixoperator{R}^{k+1}\,\frac{\velocitydvmidnext + \velocitydvmidb}{2} - \matrixoperator{P}\,\totalpdvnext =0, \\
					\matrixoperator{D}\,\left(\velocitydvmidb\right) = 0,
				\end{dcases} \label{eq::velocity_reversibility_proof_5}
			\end{equation}
			If we assume for now that $\totalpdvnext$ is known and equal to the one obtained in the forward step, it is possible to use \eqref{eq::velocity_reversibility_proof_5}  to write an expression equivalent to \eqref{eq:velocity_reversibility_proof_4} but for the backward step
			\begin{equation}
				\velocitydvmidb = \matrixoperator{B}_{-}^{-1}\matrixoperator{B}_{+}\,\velocitydvmidnext - \Delta t \,\matrixoperator{B}_{-}^{-1}\matrixoperator{P}\,\totalpdvnext.  \label{eq:velocity_reversibility_proof_6}
			\end{equation}
			Replacing \eqref{eq:velocity_reversibility_proof_4} into \eqref{eq:velocity_reversibility_proof_6} yields
			\begin{equation}
				\velocitydvmidb = \matrixoperator{B}_{-}^{-1}\matrixoperator{B}_{+}\,\left(\matrixoperator{B}_{+}^{-1}\matrixoperator{B}_{-}\,\velocitydvmidb + \Delta t \,\matrixoperator{B}_{+}^{-1}\matrixoperator{P}\,\totalpdvnext\right) + \Delta t \,\matrixoperator{B}_{-}^{-1}\matrixoperator{P}\,\totalpdvnext,  \label{eq:velocity_reversibility_proof_7}
			\end{equation}
			which can be simplified to
			\begin{equation}
				\velocitydvmidb =\velocitydvmidb +\Delta t\,\matrixoperator{B}_{-}^{-1}\matrixoperator{P}\,\left(\totalpdvnext - \totalpdvnext\right).  \label{eq:velocity_reversibility_proof_7}
			\end{equation}
			The last term on the right will be trivially equal to zero proving reversibility of the velocity time step. This term cancels out because we assumed the total pressure computed with the backward time step is equal to the one computed in the forward step. For this reason, to finalise the proof we need to show that indeed the forward step total pressure is the solution for the pressure in the backward step. We start with a modified version of \eqref{eq::velocity_reversibility_proof_2}
			\begin{equation}
				\begin{dcases}
					\matrixoperator{M} \frac{\velocitydvmidnext - \velocitydvmidb}{\Delta t} + \matrixoperator{R}^{k+1}\,\frac{\velocitydvmidnext + \velocitydvmidb}{2} - \matrixoperator{P}\,\totalpdvnext =0, \\
					\matrixoperator{D}\,\left(\velocitydvmidnext + \velocitydvmidb\right) = 0.
				\end{dcases} \label{eq::velocity_reversibility_proof_8}
			\end{equation}
			This system of equations is equivalent to $\eqref{eq::velocity_reversibility_proof_2}$ because the divergence of the discrete velocity is zero, therefore $\matrixoperator{D}\,\velocitydmidb = 0$. Rearranging gives
			\begin{equation}
				\begin{dcases}
					\left(\matrixoperator{M} + \frac{\Delta t}{2}\matrixoperator{R}^{k+1}\right)\,\velocitydvmidnext  - \left(\matrixoperator{M} - \frac{\Delta t}{2}\matrixoperator{R}^{k+1}\right)\,\velocitydvmidb - \Delta t\,\matrixoperator{P}\,\totalpdvnext =0, \\
					\matrixoperator{D}\,\left(\velocitydvmidnext + \velocitydvmidb\right) = 0.
				\end{dcases} \label{eq::velocity_reversibility_proof_8}
			\end{equation}
			In the same way the backwards step results in the following system of equations
			\begin{equation}
				\begin{dcases}
					-\left(\matrixoperator{M} - \frac{\Delta t}{2}\matrixoperator{R}^{k+1}\right)\,\velocitydvmidb  + \left(\matrixoperator{M} + \frac{\Delta t}{2}\matrixoperator{R}^{k+1}\right)\,\velocitydvmidnext - \Delta t\,\matrixoperator{P}\,\totalpdvnext =0, \\
					\matrixoperator{D}\,\left(\velocitydvmidnext + \velocitydvmidb\right) = 0.
				\end{dcases} \label{eq::velocity_reversibility_proof_9}
			\end{equation}
			The system \eqref{eq::velocity_reversibility_proof_8} is the same as \eqref{eq::velocity_reversibility_proof_9}, therefore the pressure is necessarily the same in both cases. Thus proving the reversibility of the velocity time step.
				
\section{Numerical test cases}  \label{sec::numerical_test_cases}
	In order to verify the numerical properties of the proposed method we apply it to the solution of two well known test cases: (i) Taylor-Green vortex and (ii) inviscid shear layer roll-up. With the first test the convergence properties of the method will be assessed and with the second the time reversibility and mass, energy, enstrophy, and vorticity conservation will be verified.
		
	\subsection{Order of accuracy study: Taylor-Green vortex}\label{subsec::taylor_green_vortex}
		In order to verify the accuracy of the proposed numerical method we compare the numerical results against an analytical solution of the Navier-Stokes equations. A suitable analytical solution is the Taylor-Green vortex, given by:
		\begin{equation}
			\begin{dcases}
				u_{x} (x,y,t) = -\sin(\pi x)\cos(\pi y)\,e^{-2\pi^{2} \nu t}, \\
				u_{y} (x,y,t) = \cos(\pi x)\sin(\pi y)\,e^{-2\pi^{2} \nu t}, \\
				p(x,y,t) = \frac{1}{4}\,(\cos(2\pi x)+ \cos(2\pi y))\,e^{4\pi^{2} \nu t}, \\
				\omega (x,y,t) = -2\pi\sin(\pi x)\sin(\pi y)\,e^{-2\pi^{2} \nu t}.
			\end{dcases} \label{eq:taylor_green_exact_solution}
		\end{equation}
		The solution is defined on the domain $\Omega = [0,2]\times[0,2]$, with periodic boundary conditions. The initial condition for both the velocity $\velocitydstart$ and vorticity $\vorticitydstart$ are given by the exact solution \eqref{eq:taylor_green_exact_solution}. The kinematic viscosity is set to $\nu = 0.01$. For this study we consider the evolution of the solution from $t=0$ to $t=1$.
		
\begin{figure}[!ht]
			\centering
			\includegraphics[width=0.45\textwidth]{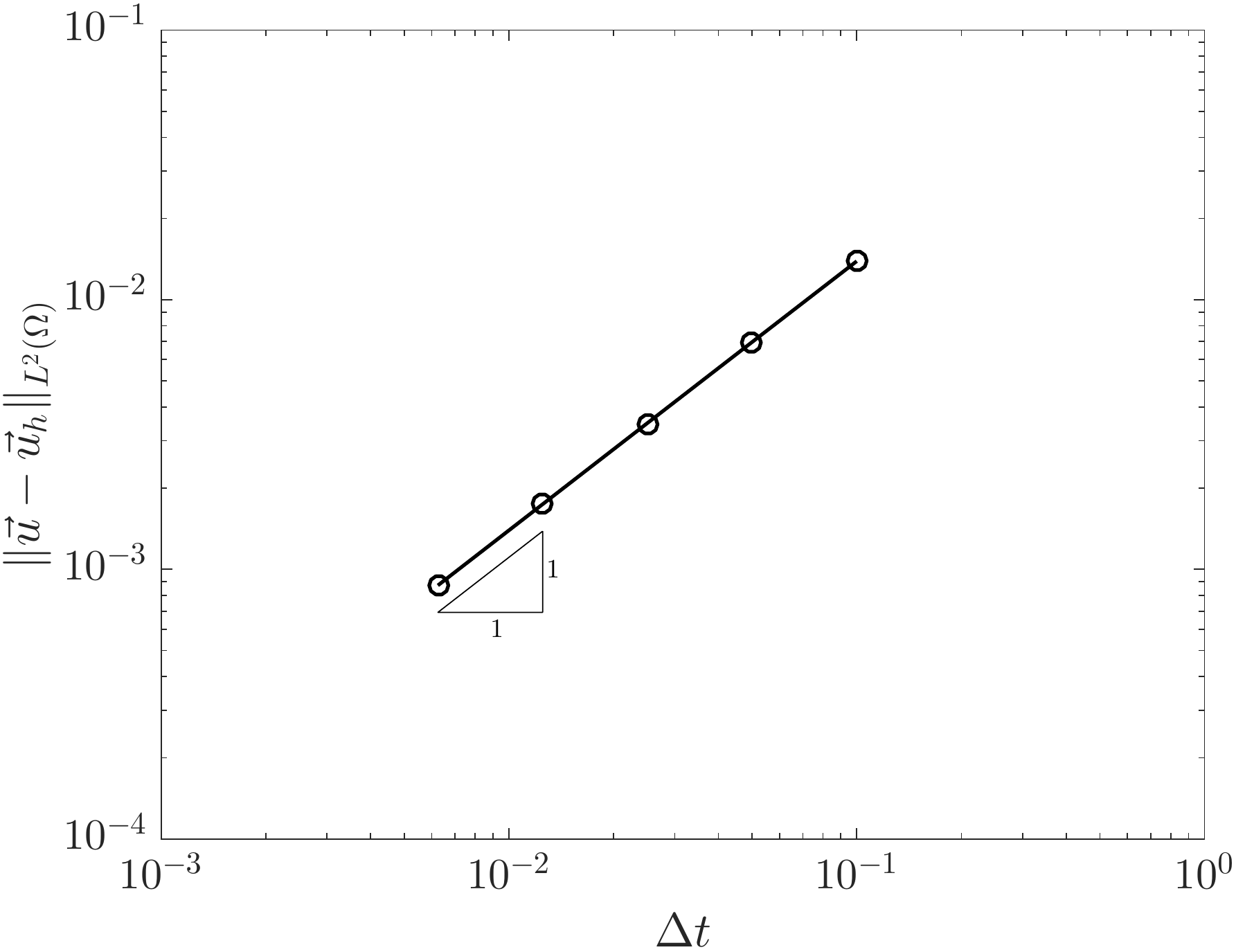}
			\includegraphics[width=0.45\textwidth]{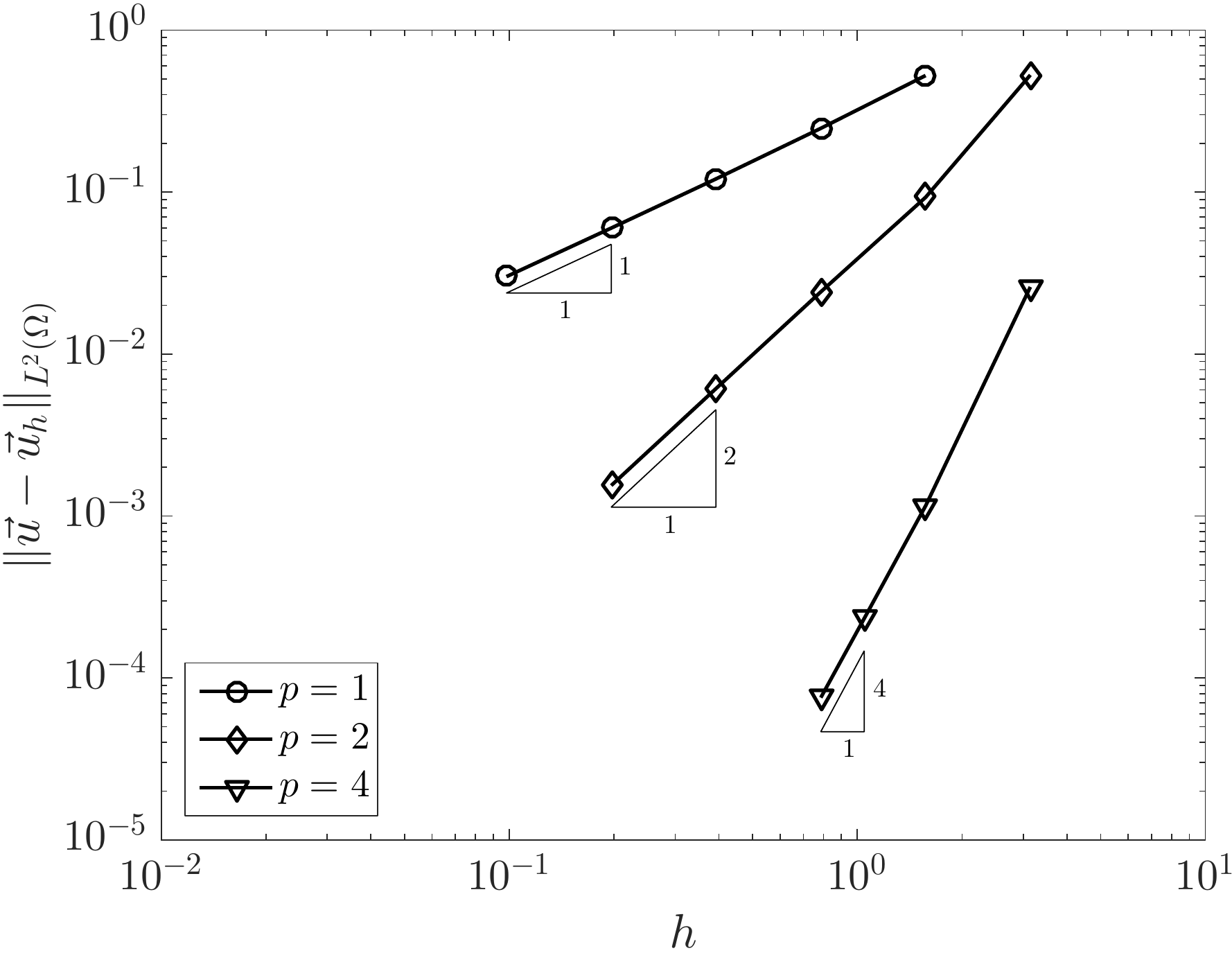}
			\caption{Convergence plots of velocity error, $\|\velocity-\velocityd\|_{L^{2}(\Omega)}$. Left: error in velocity as function of time step size $\Delta t$. Right: error in velocity as function of mesh size, $h$, for basis functions of polynomial degree $p=1,2,4$.}
			\label{fig:convergence_plots}
		\end{figure}

		The first study focusses on the time convergence. For this we have used 1024 triangular elements and polynomial basis functions of degree $p=4$ and time steps equal to $\Delta t = 1, \frac{1}{2}, \frac{1}{4}, \frac{1}{8}, \frac{1}{16}$. As can be observed in \figref{fig:convergence_plots} left, this method achieves a first order convergence rate, as opposed to a second order convergence rate of an implicit formulation, see \cite{Sanderse2013}.
		
		Regarding the spatial convergence we have tested the convergence rate of discretizations with basis functions of different polynomial degree, $p=1,2,4$. In order not to pollute the spatial convergence rate with the temporal integration error, we have used different time steps: $\Delta t = 2.5\times 10^{-2}$ for $p=1$, $\Delta t = 1.0\times 10^{-3}$ for $p=2$ and $\Delta t = 1.0\times 10^{-4}$ for $p=4$. As can be seen in \figref{fig:convergence_plots}, this method has a convergence rate of $p$-th order for basis functions of polynomial degree $p$.
		
	\subsection{Time reversibility and mass, energy, enstrophy, and vorticity conservation: inviscid shear layer roll-up} \label{subsec::inviscid_shear_layer_roll_up}
		The second group of tests focusses on the conservation properties of the proposed numerical method. For this set of tests we consider inviscid flow $\nu=0$. We consider here the simulation of the roll-up of a shear layer, see e.g. \cite{Minion1997,Knikker2009,Hokpunna2010,Sanderse2013}. This solution is particularly challenging because during the evolution large vorticity gradients develop. Several methods are known to \emph{blow up}, e.g. \cite{Minion1997}. We consider on the periodic domain $\Omega = [0,2\pi]\times[0,2\pi]$ the following initial conditions for velocity $\velocity$ 
		\begin{equation}
			u_{x}(x,y) = 
				\begin{dcases}
					\tanh\left(\frac{y-\frac{\pi}{2}}{\delta}\right), & y\leq\pi, \\
					\tanh\left(\frac{\frac{3\pi}{2}-y}{\delta}\right), & y>\pi,
				\end{dcases} 
			\qquad\qquad\qquad
			u_{y}(x,y) = \epsilon\sin(x),
		\end{equation}
		and vorticity $\vorticity$
		\begin{equation}
			\omega(x,y) = 
				\begin{dcases}
					\frac{1}{\delta}\,\text{sech}^2\left(\frac{y-\frac{\pi }{2}}{\delta }\right), & y\leq\pi, \\
					-\frac{1}{\delta}\,\text{sech}^2\left(\frac{\frac{3 \pi }{2}-y}{\delta }\right), & y>\pi, 
				\end{dcases}
		\end{equation}
		with $\delta = \frac{\pi}{15}$ and $\epsilon=0.05$ as in \cite{Knikker2009,Sanderse2013}.
		
\begin{figure}[!ht]
			\centering
			\includegraphics[width=0.4\textwidth]{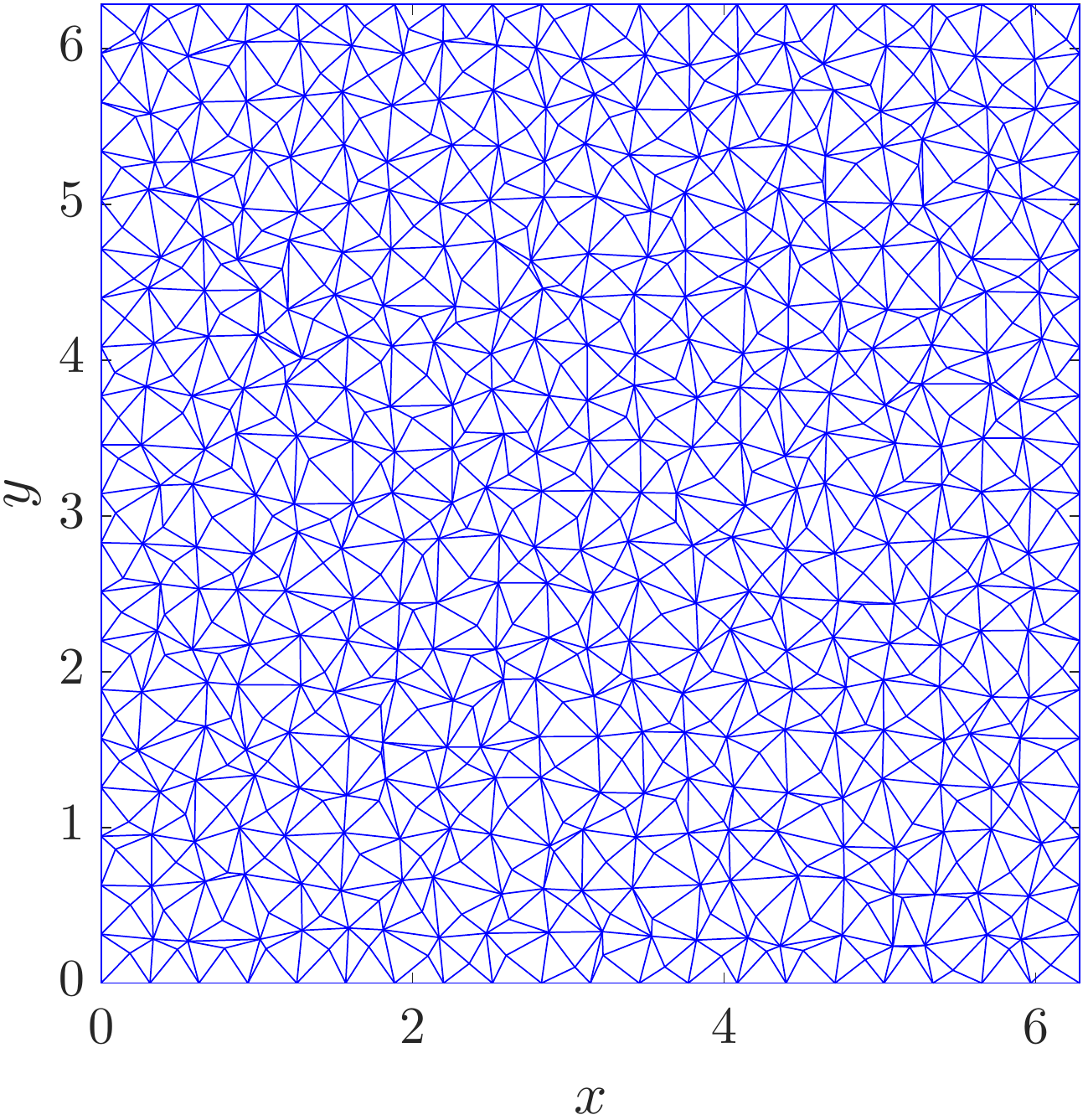}
			\caption{Example mesh with 1600 triangular elements.}
			\label{fig:mesh}
		\end{figure}

\begin{figure}[!htb]
			\centering
			\includegraphics[width=0.45\textwidth]{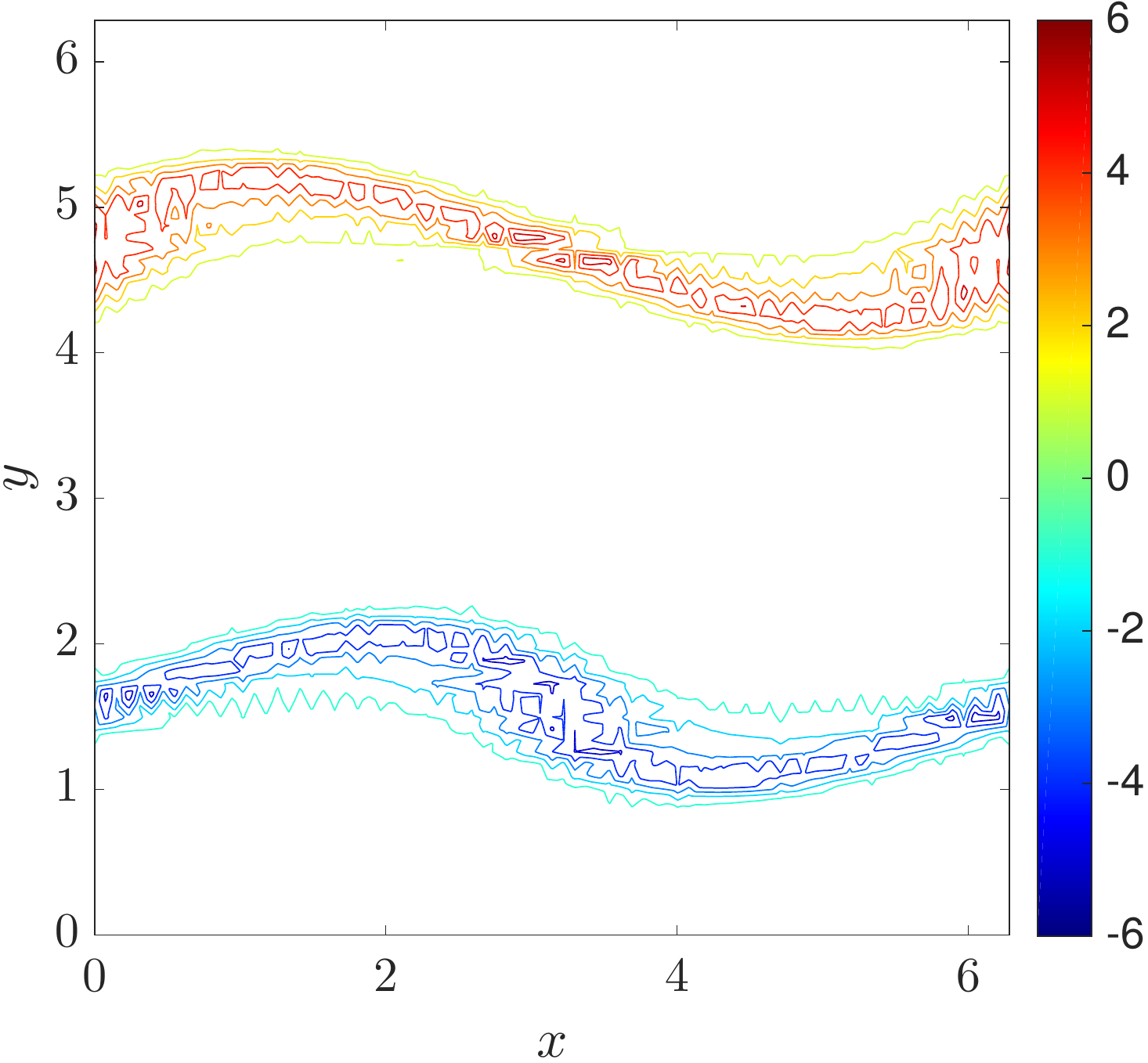}
			\includegraphics[width=0.45\textwidth]{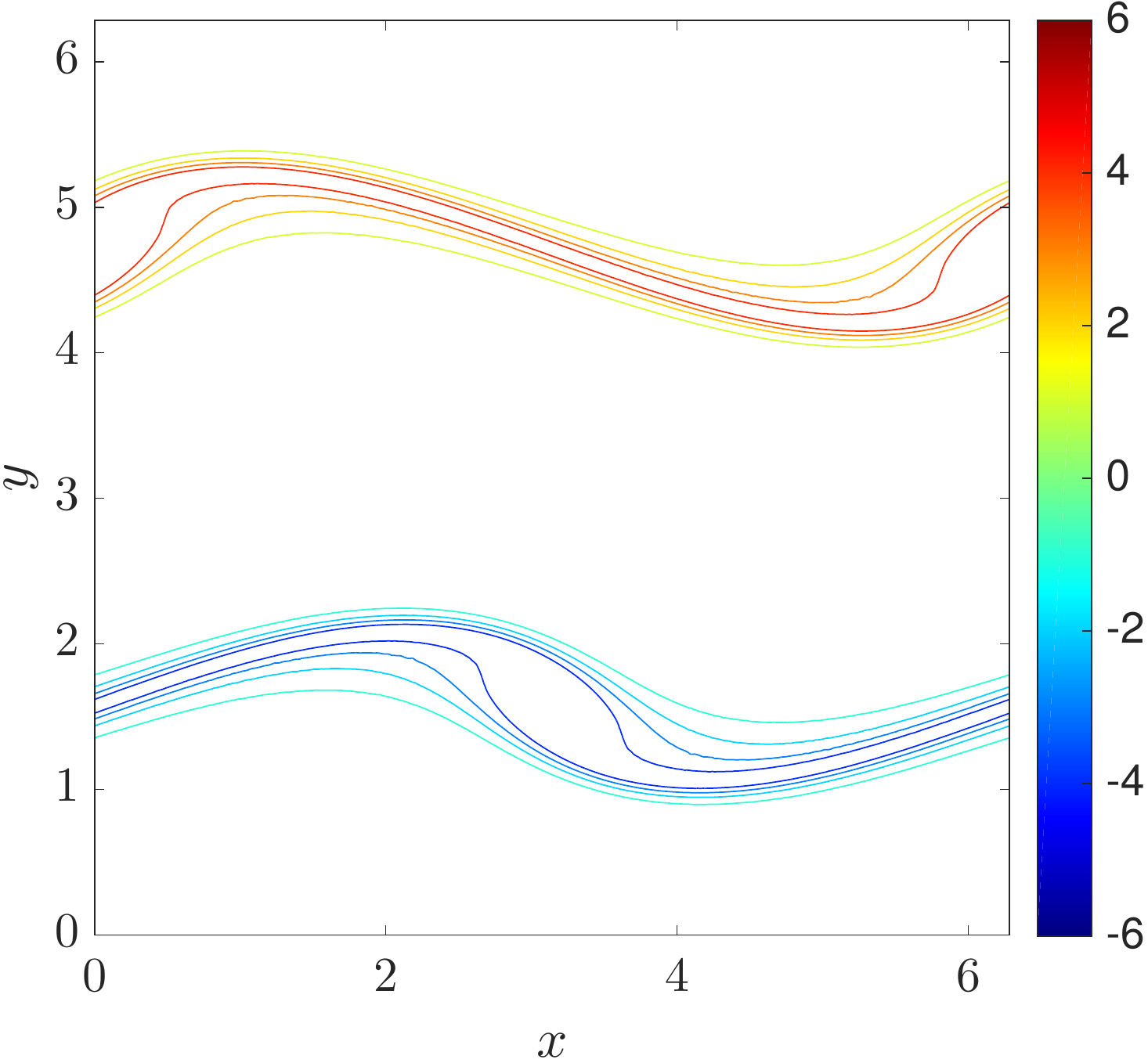}
			\caption{Vorticity field of inviscid shear layer roll-up at $t=4$ obtained with 6400 triangular finite elements, contour lines ($\pm 1, \pm 2, \pm 3, \pm 4, \pm 5, \pm 6$) as in \cite{Knikker2009}. Left: solution with basis functions of polynomial degree $p=1$ and time step $\Delta t = 0.1$. Right: solution with basis functions of polynomial degree $p=4$ and time step $\Delta t = 0.01$.}
			\label{fig:vorticity_evolution_comparison_plots_contour_1}
		\end{figure}

\begin{figure}[!htb]
			\centering
			\includegraphics[width=0.45\textwidth]{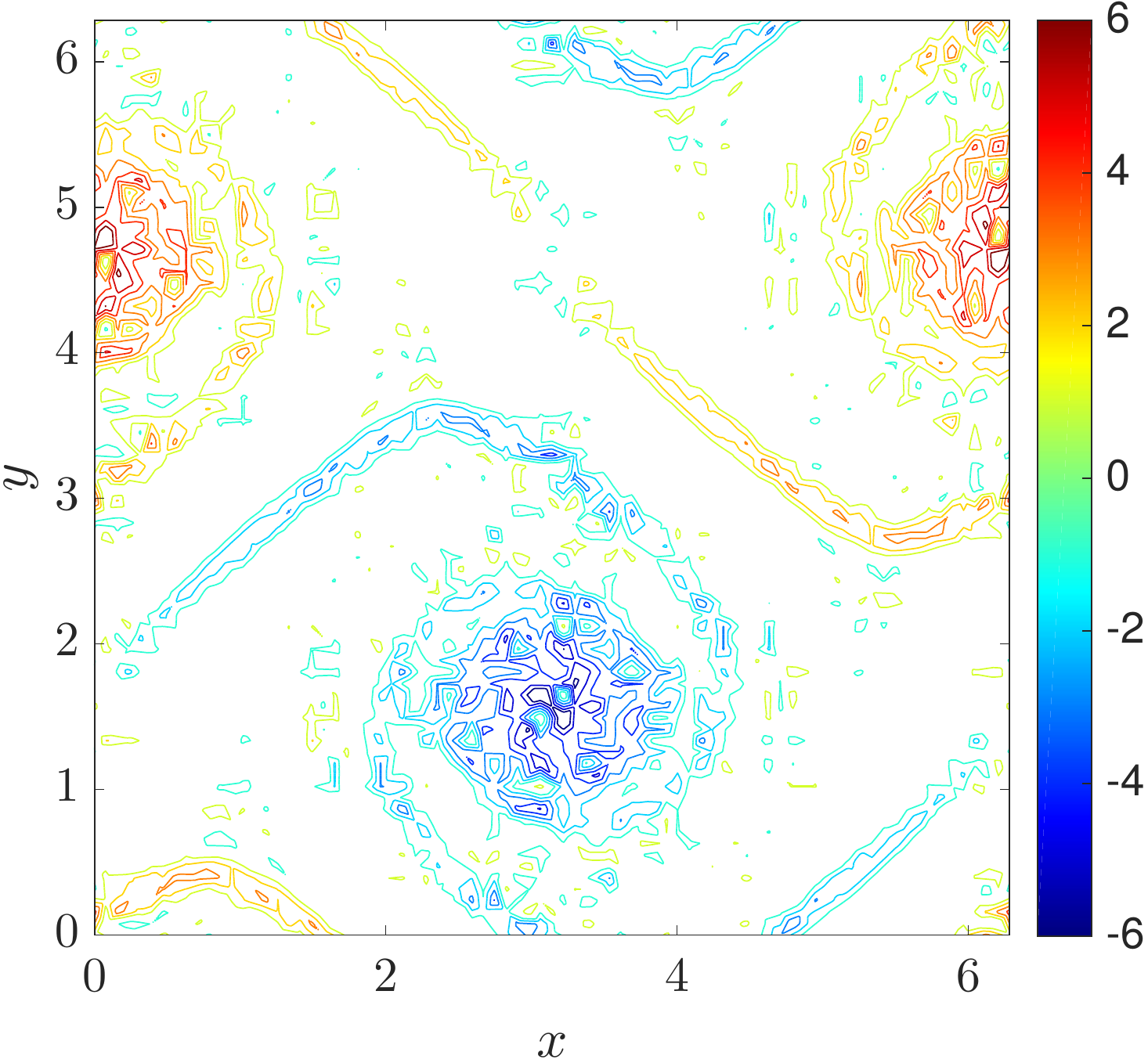}
			\includegraphics[width=0.45\textwidth]{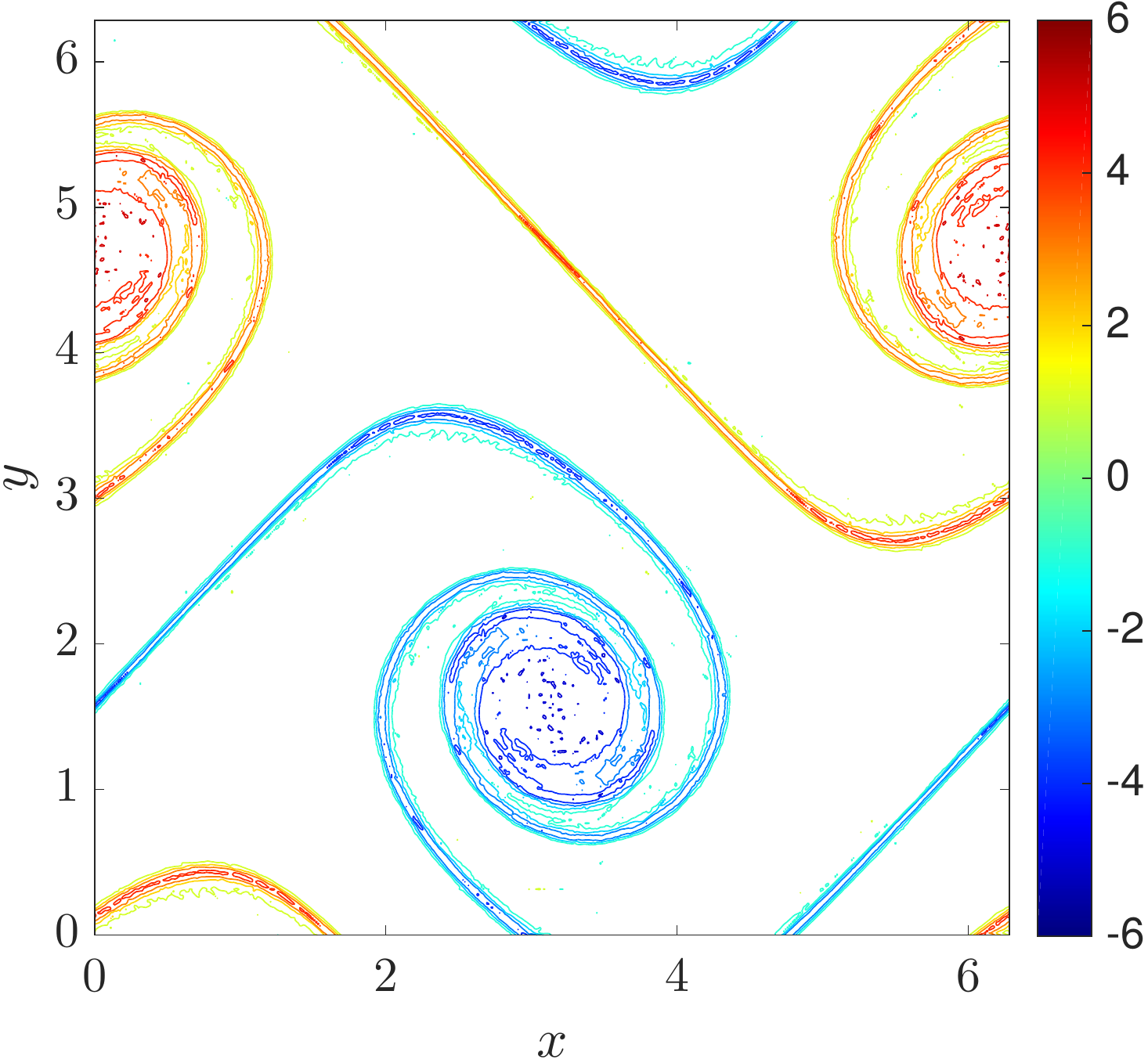}
			\caption{Vorticity field of inviscid shear layer roll-up at $t=8$ obtained with 6400 triangular finite elements, contour lines ($\pm 1, \pm 2, \pm 3, \pm 4, \pm 5, \pm 6$) as in \cite{Knikker2009}. Left: solution with basis functions of polynomial degree $p=1$ and time step $\Delta t = 0.1$s. Right: solution with basis functions of polynomial degree $p=4$ and time step $\Delta t = 0.01$s.}
			\label{fig:vorticity_evolution_comparison_plots_contour_2}
		\end{figure}

		The small perturbation $\epsilon$ in the $y$-component of the velocity field will trigger the roll-up of the shear layer. We show the contour lines of vorticity at $t=4$ and $t=8$, \figref{fig:vorticity_evolution_comparison_plots_contour_1} and \figref{fig:vorticity_evolution_comparison_plots_contour_2} respectively. The plots on the left side of \figref{fig:vorticity_evolution_comparison_plots_contour_1} and \figref{fig:vorticity_evolution_comparison_plots_contour_2} correspond to a spatial discretization of 6400 triangular elements of polynomial degree $p=1$ and time step size $\Delta t = 0.1$. On the right we present the results for a spatial discretization of 6400 triangular elements of polynomial degree $p=4$ and time step size $\Delta t = 0.01$. An example mesh with 1600 elements is shown in \figref{fig:mesh}.
		
		Although it is possible to observe oscillations in the vorticity plots, especially for $p=1$, none of them \emph{blows up}, as is reported in \cite{Minion1997}. These oscillations necessarily appear whenever the flow generates structures at a scale smaller than the resolution of the mesh. Since energy and enstrophy are conserved, there is no possibility for the numerical simulation to dissipate the small scale information. This is a feature observed in all conserving inviscid solvers without a small scale model.
		
		
\begin{figure}[!hb]
			\centering
			\includegraphics[width=0.4\textwidth]{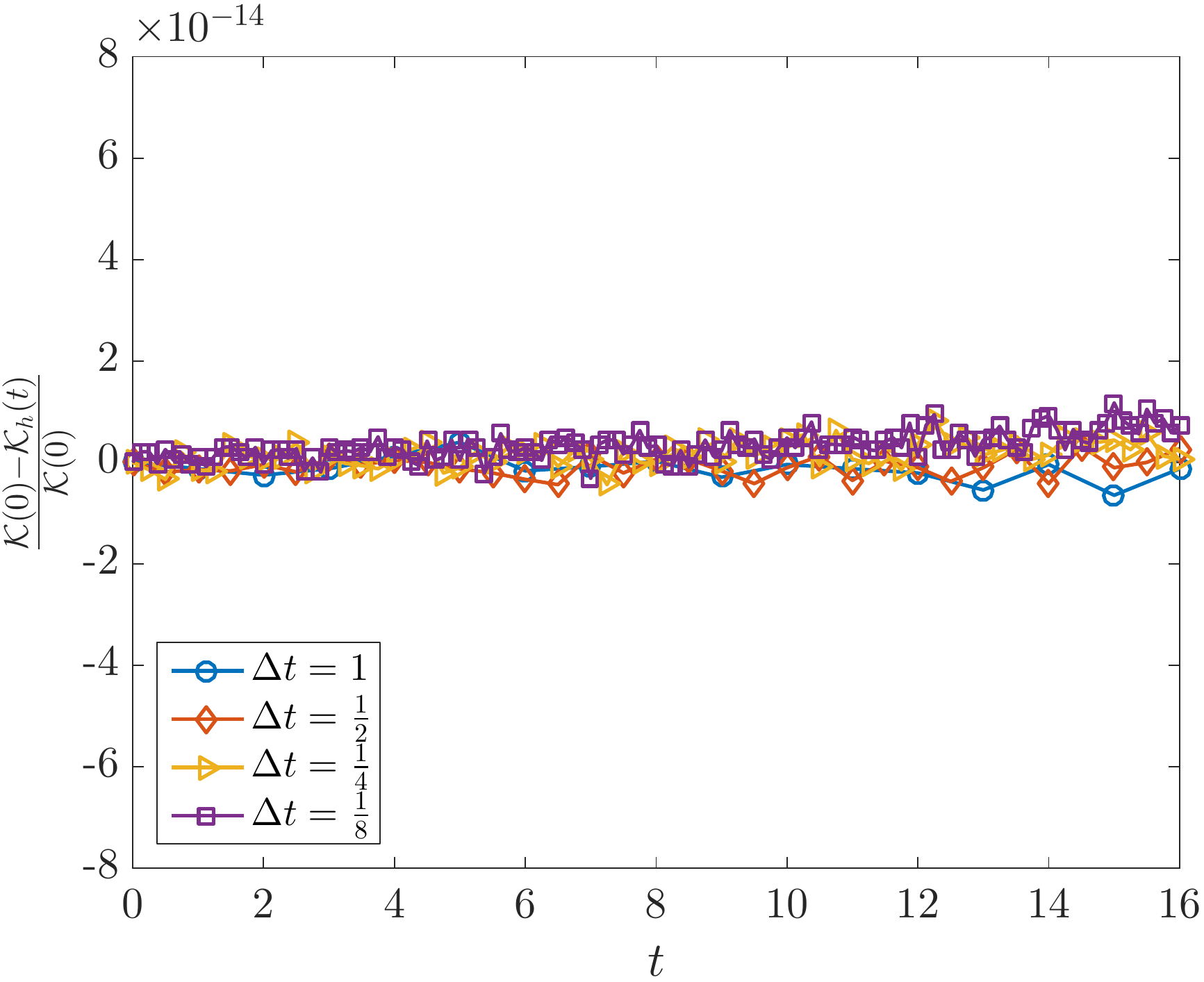}
			\includegraphics[width=0.4\textwidth]{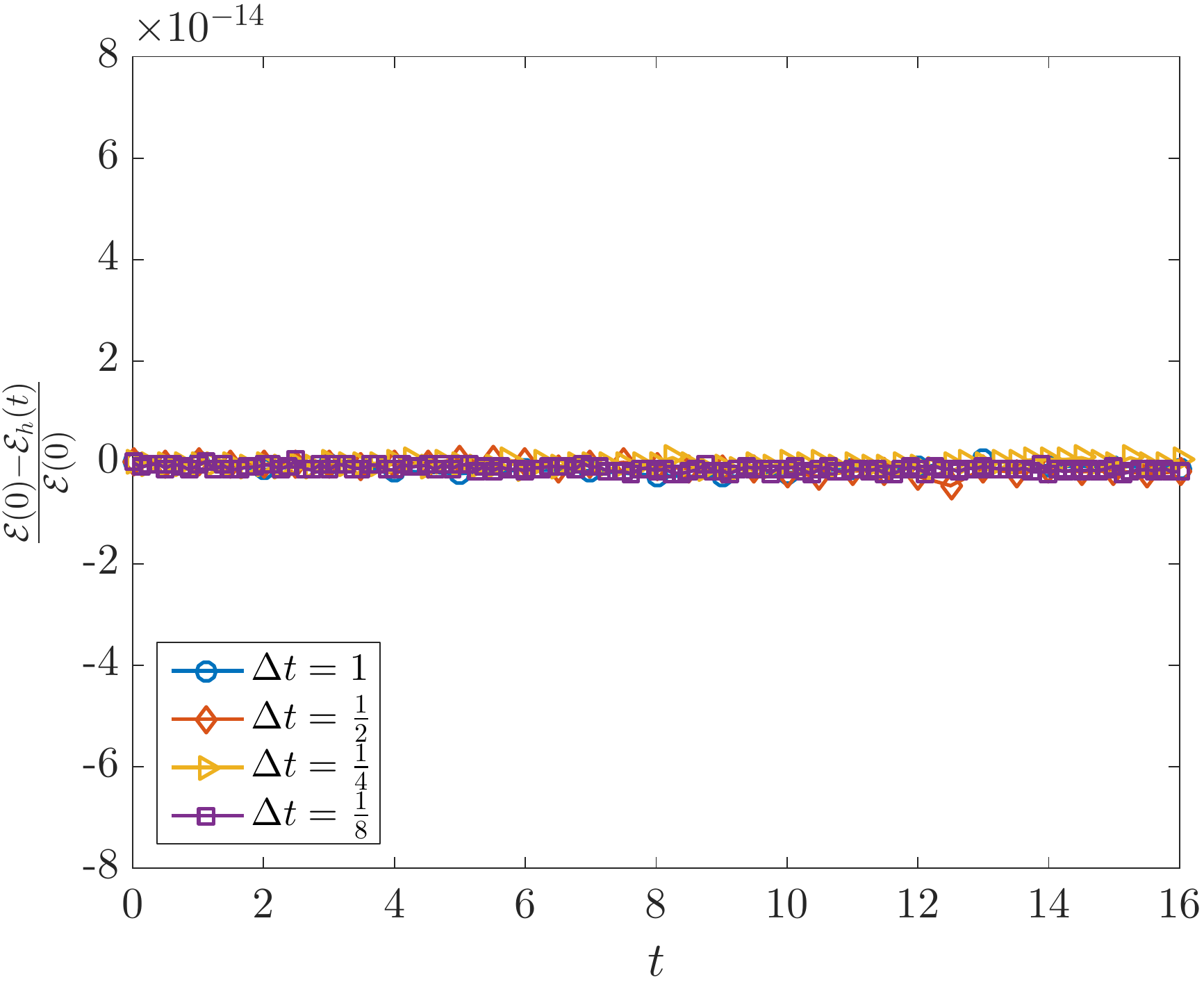}
			\caption{Simulation of shear-layer problem with 6400 triangular finite elements, basis functions of polynomial degree $p=1$, and different time steps $\Delta t = 1,\frac{1}{2}, \frac{1}{4},\frac{1}{8}$. Left: kinetic energy error with respect to initial kinetic energy, $\frac{\energy(0)-\energyd(t)}{\energy(0)}$. Right: enstrophy error with respect to initial enstrophy, $\frac{\enstrophy(0)-\enstrophyd(t)}{\enstrophy(0)}$.}
			\label{fig:energy_enstrophy_conservation_plots}
		\end{figure}

\begin{figure}[!htb]
			\centering
			\includegraphics[width=0.4\textwidth]{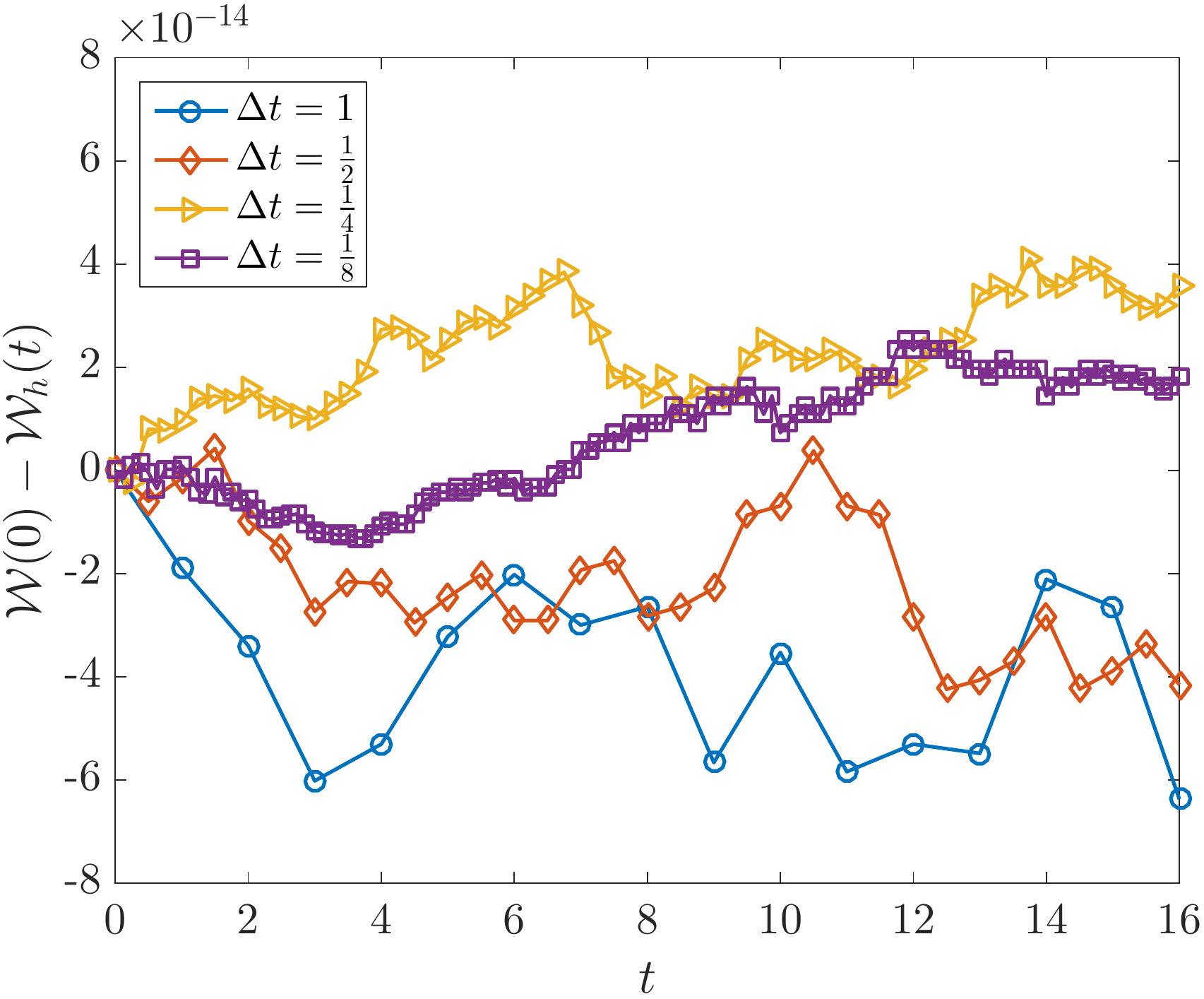}
			\includegraphics[width=0.4\textwidth]{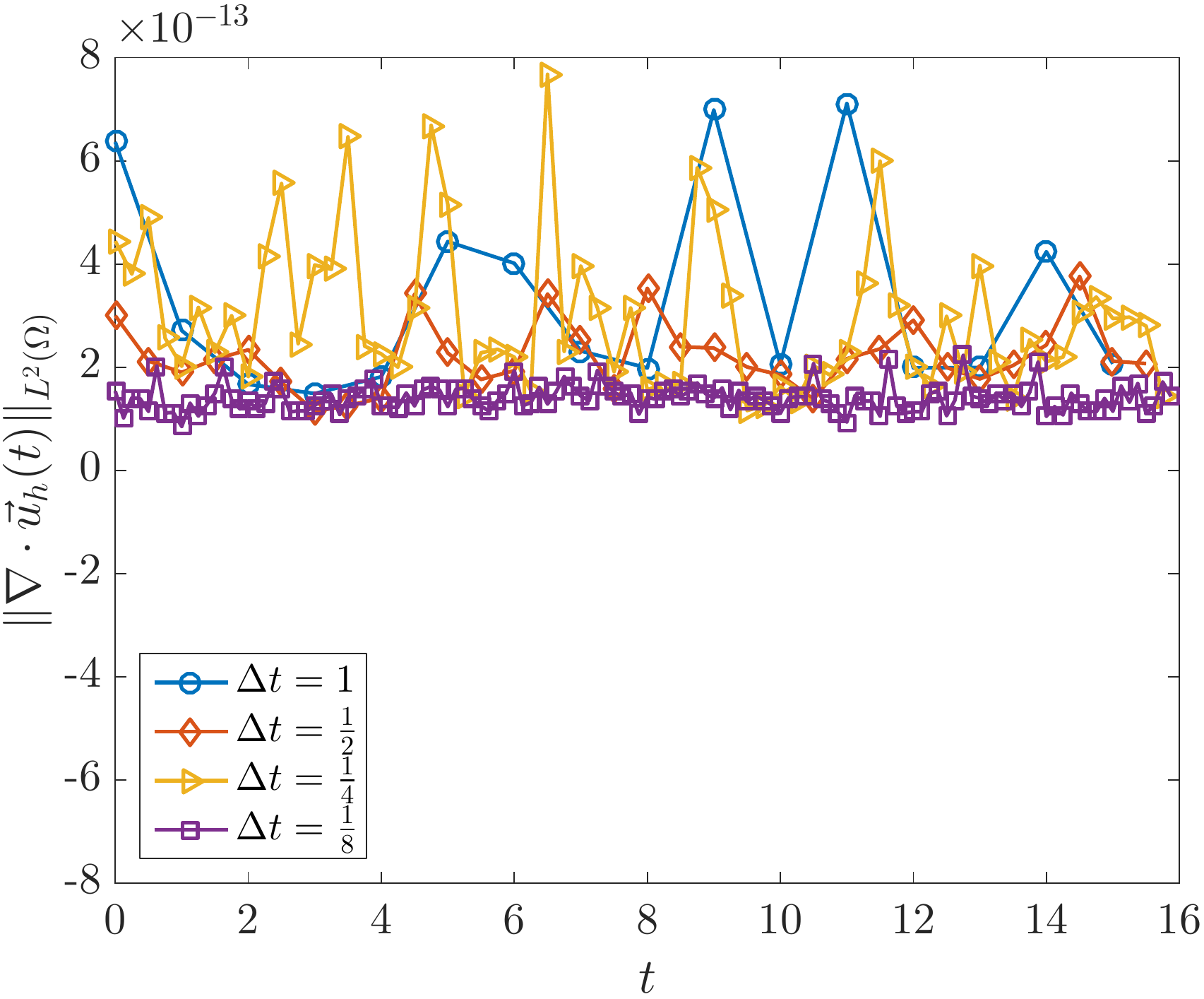}
			\caption{Simulation of shear-layer problem with 6400 triangular finite elements, basis functions of polynomial degree $p=1$, and different time steps $\Delta t = 1,\frac{1}{2}, \frac{1}{4},\frac{1}{8}$. Left: total vorticity error with respect to initial total vorticity, $\vorticityt(0)-\vorticitytd(t)$. Right: evolution of the divergence of the velocity field.}
			\label{fig:vorticity_divergence_conservation_plots}
		\end{figure}

		In order to assess the conservation properties of the proposed method we computed the evolution of the shear-layer problem from $t=0$ to $t=16$, using different time steps $\Delta t = 1,\frac{1}{2}, \frac{1}{4},\frac{1}{8}$. We have used a coarse grid with  6400 triangular finite elements and basis functions of polynomial degree $p=1$ for the spatial discretization. In \figref{fig:energy_enstrophy_conservation_plots} left we show the evolution of the kinetic energy error with respect to the initial kinetic energy,  $\frac{\energy(0)-\energyd(t)}{\energy(0)}$. This result confirms conservation of kinetic energy, as proven in \secref{sec::kinetic_energy_conservation}. A similar study was performed for enstrophy in \figref{fig:energy_enstrophy_conservation_plots} right. The magnitude of the error observed is of machine precision, verifying conservation of enstrophy, as proven in \secref{sec::enstrophy_conservation}. The evolution of the error of total vorticity is shown in \figref{fig:vorticity_divergence_conservation_plots} left. The error of vorticity shows again an error of the order of machine precision. In \figref{fig:vorticity_divergence_conservation_plots} right we show the divergence of the velocity field. As can be seen, this discretization results in a velocity field that is divergence free.
		
\begin{figure}[!hb]
			\centering
			\includegraphics[width=0.4\textwidth]{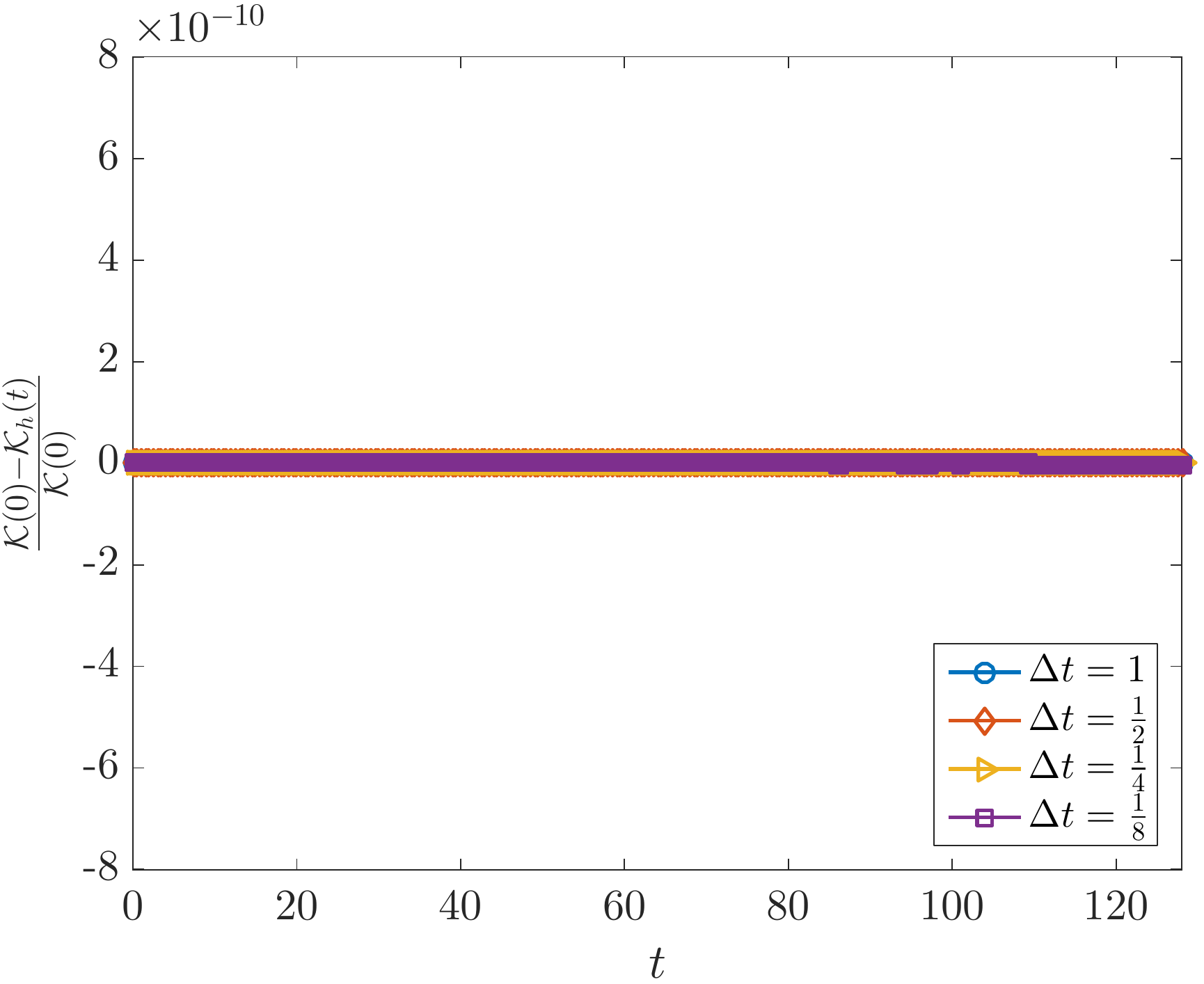}
			\includegraphics[width=0.4\textwidth]{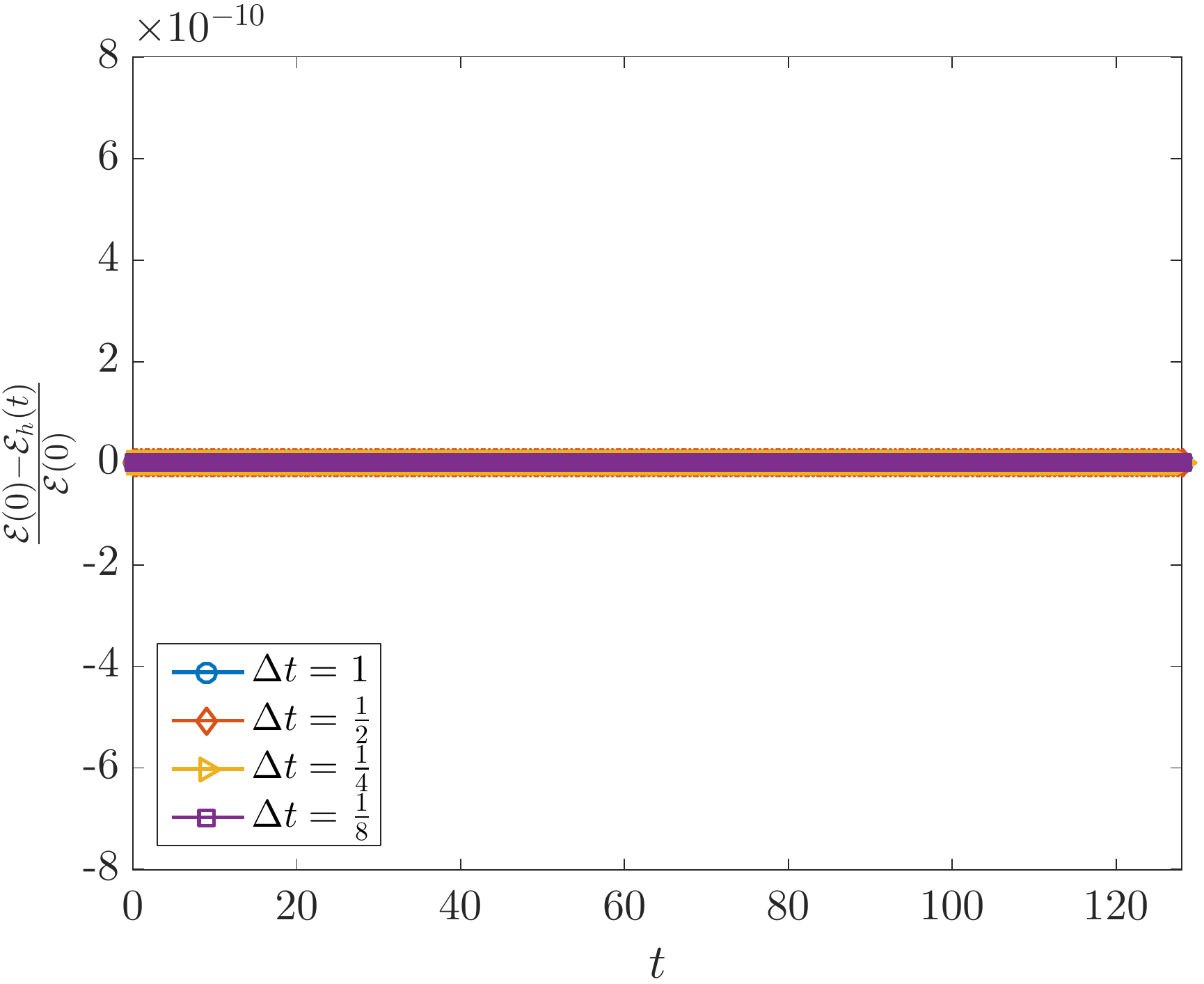}
			\caption{Simulation of shear-layer problem with 100 triangular finite elements, basis functions of polynomial degree $p=4$, and different time steps $\Delta t = 1,\frac{1}{2}, \frac{1}{4},\frac{1}{8}$. Left: kinetic energy error with respect to initial kinetic energy, $\frac{\energy(0)-\energyd(t)}{\energy(0)}$. Right: enstrophy error with respect to initial enstrophy, $\frac{\enstrophy(0)-\enstrophyd(t)}{\enstrophy(0)}$.}
			\label{fig:energy_enstrophy_conservation_plots_long_time}
		\end{figure}

\begin{figure}[!hb]
			\centering
			\includegraphics[width=0.4\textwidth]{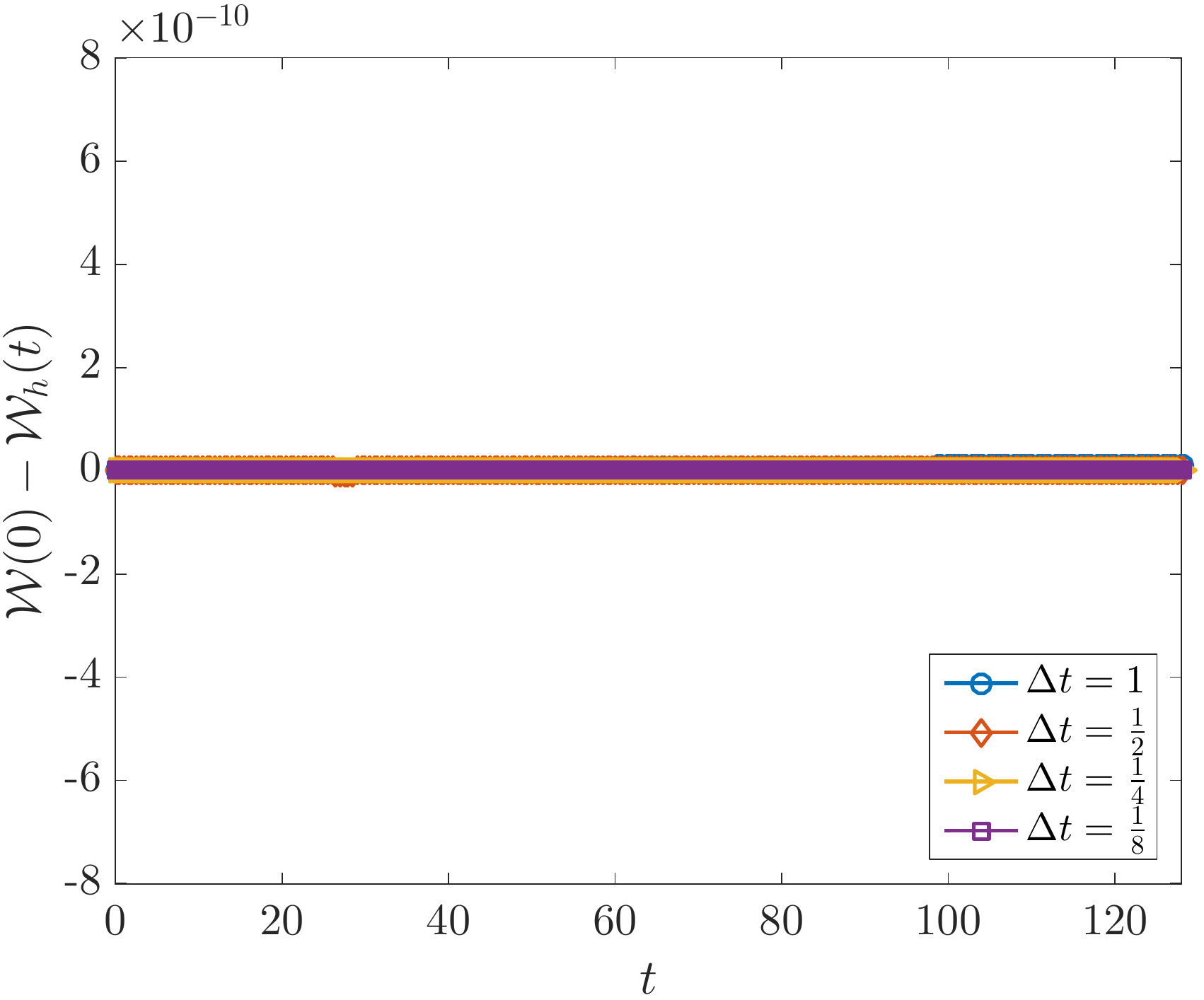}
			\includegraphics[width=0.4\textwidth]{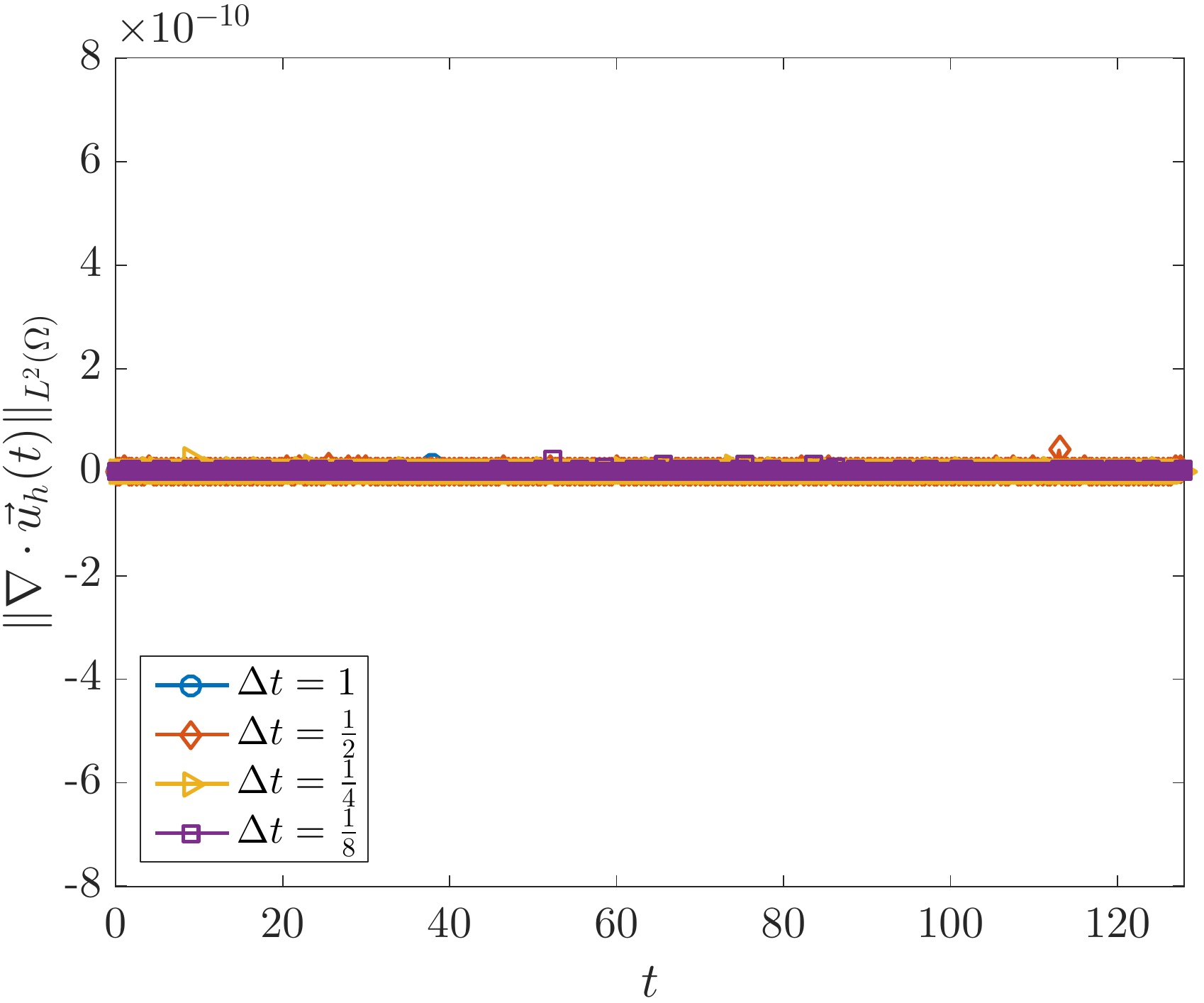}
			\caption{Simulation of shear-layer problem with 100 triangular finite elements, basis functions of polynomial degree $p=4$, and different time steps $\Delta t = 1,\frac{1}{2}, \frac{1}{4},\frac{1}{8}$. Left: total vorticity error with respect to initial total vorticity, $\vorticityt(0)-\vorticitytd(t)$. Right: evolution of the divergence of the velocity field.}
			\label{fig:vorticity_divergence_conservation_plots_long_time}
		\end{figure}

		In order to assess the conservation properties of the proposed method for long time simulations we computed the evolution of the shear-layer problem from $t=0$ to $t=128$, using different time steps $\Delta t = 1,\frac{1}{2}, \frac{1}{4},\frac{1}{8}$. We have used a very coarse grid with  100 triangular finite elements and basis functions of polynomial degree $p=4$ for the spatial discretization. In \figref{fig:energy_enstrophy_conservation_plots_long_time} left we show the evolution of the kinetic energy error with respect to the initial kinetic energy,  $\frac{\energy(0)-\energyd(t)}{\energy(0)}$. This result confirms conservation of kinetic energy, as proven in \secref{sec::kinetic_energy_conservation}. The same study was performed for enstrophy in \figref{fig:energy_enstrophy_conservation_plots_long_time} right. The error observed is very small, verifying conservation of enstrophy, as proven in \secref{sec::enstrophy_conservation}. The evolution of error for total vorticity is shown in \figref{fig:vorticity_divergence_conservation_plots_long_time} left. The error of total vorticity shows that this method is capable of conserving the total vorticity on long time simulations. In \figref{fig:vorticity_divergence_conservation_plots_long_time} right we show the divergence of the velocity field. As can be seen, this discretization results in a velocity field that is divergence free. Similar results are obtained for $p=1$.
		
\begin{figure}[!htb]
			\centering
			\begin{minipage}{\textwidth}\footnotesize
				\centering
				\subsubfloat{\includegraphics[width=0.45\textwidth]{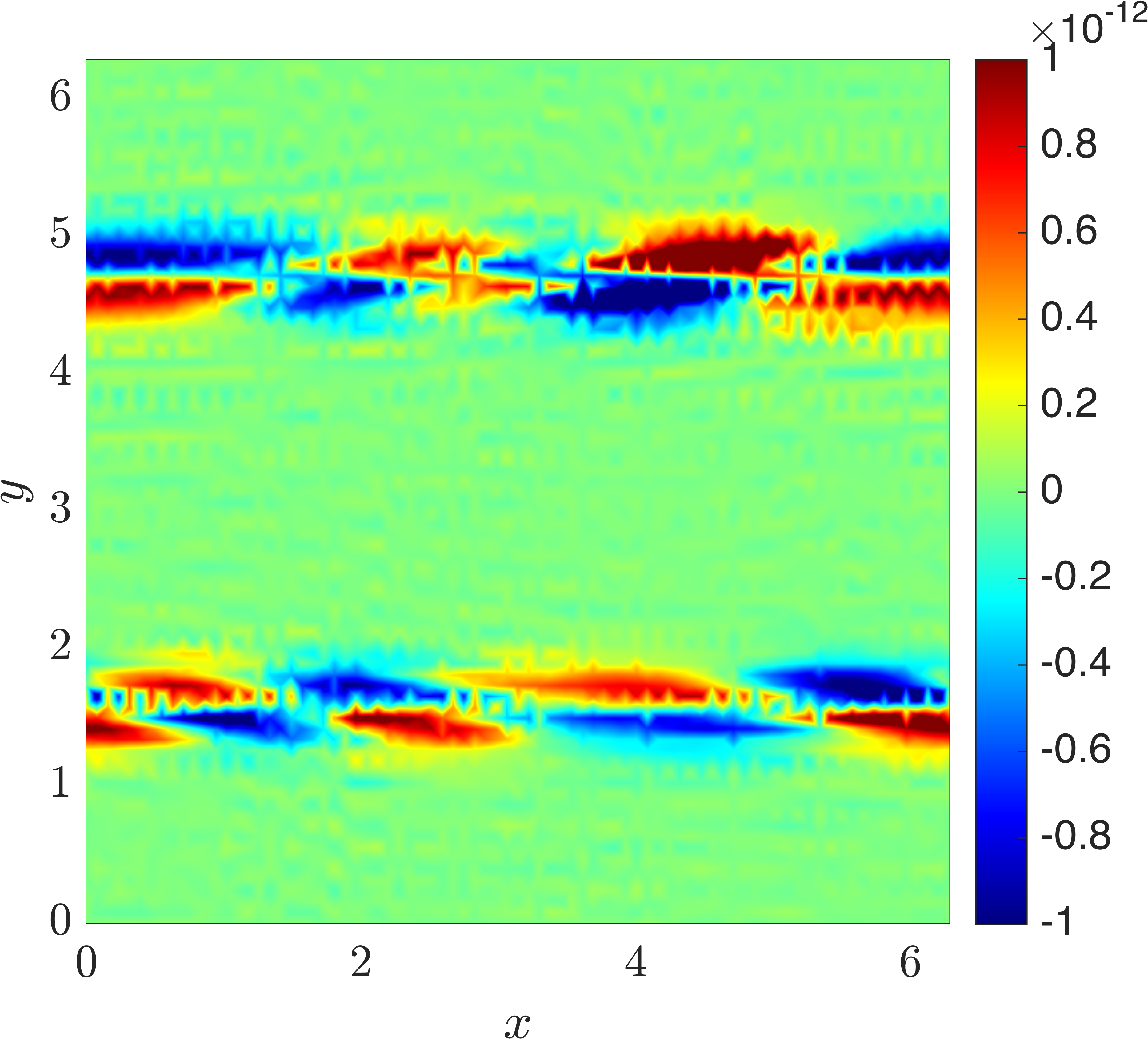}}{a) $\Delta t = \frac{1}{8}$.}
				\subsubfloat{\includegraphics[width=0.45\textwidth]{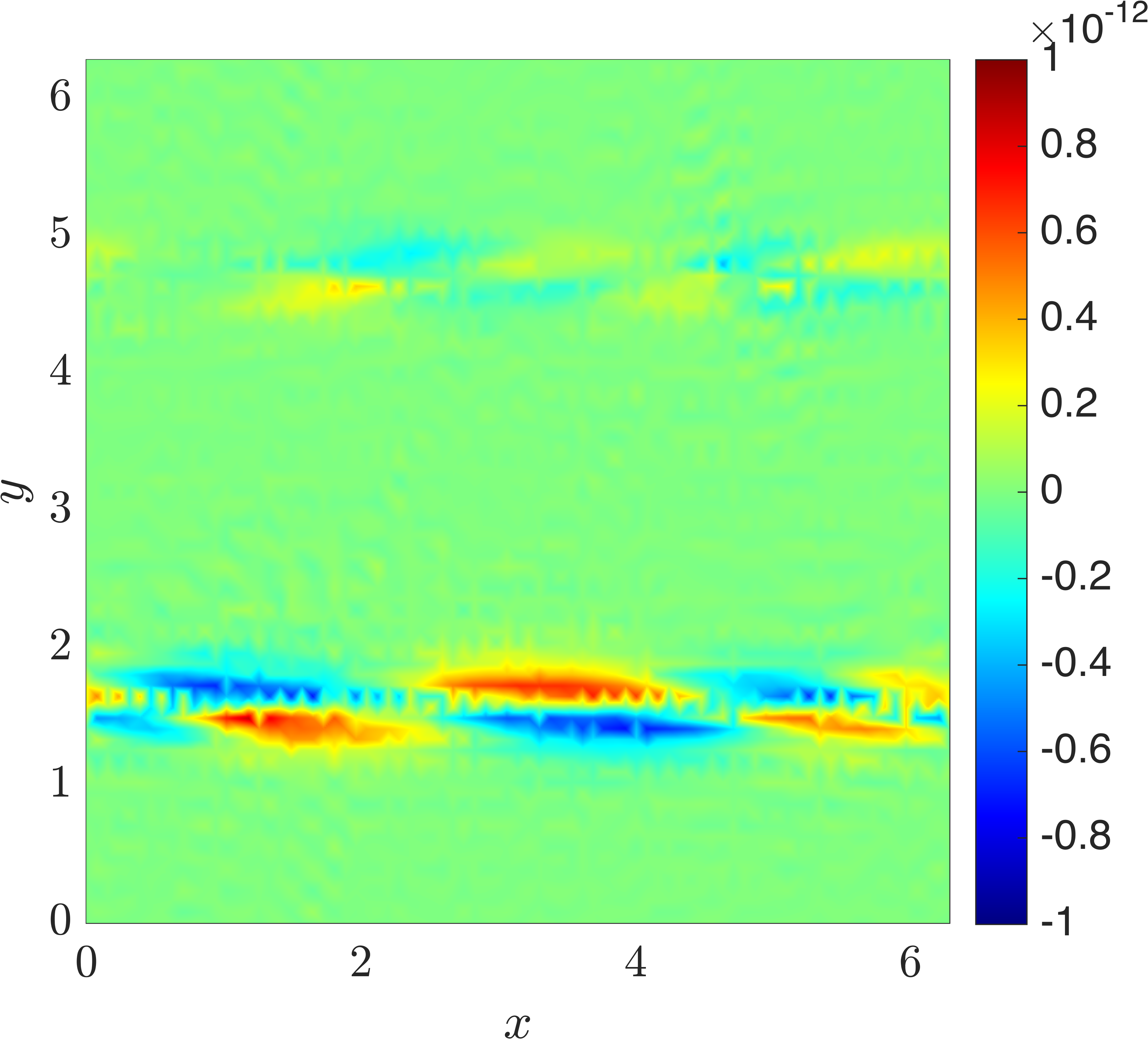}}{b) $\Delta t = \frac{1}{4}$.}
			\end{minipage}
			
			\begin{minipage}{\textwidth}\footnotesize
				\centering
				\subsubfloat{\includegraphics[width=0.45\textwidth]{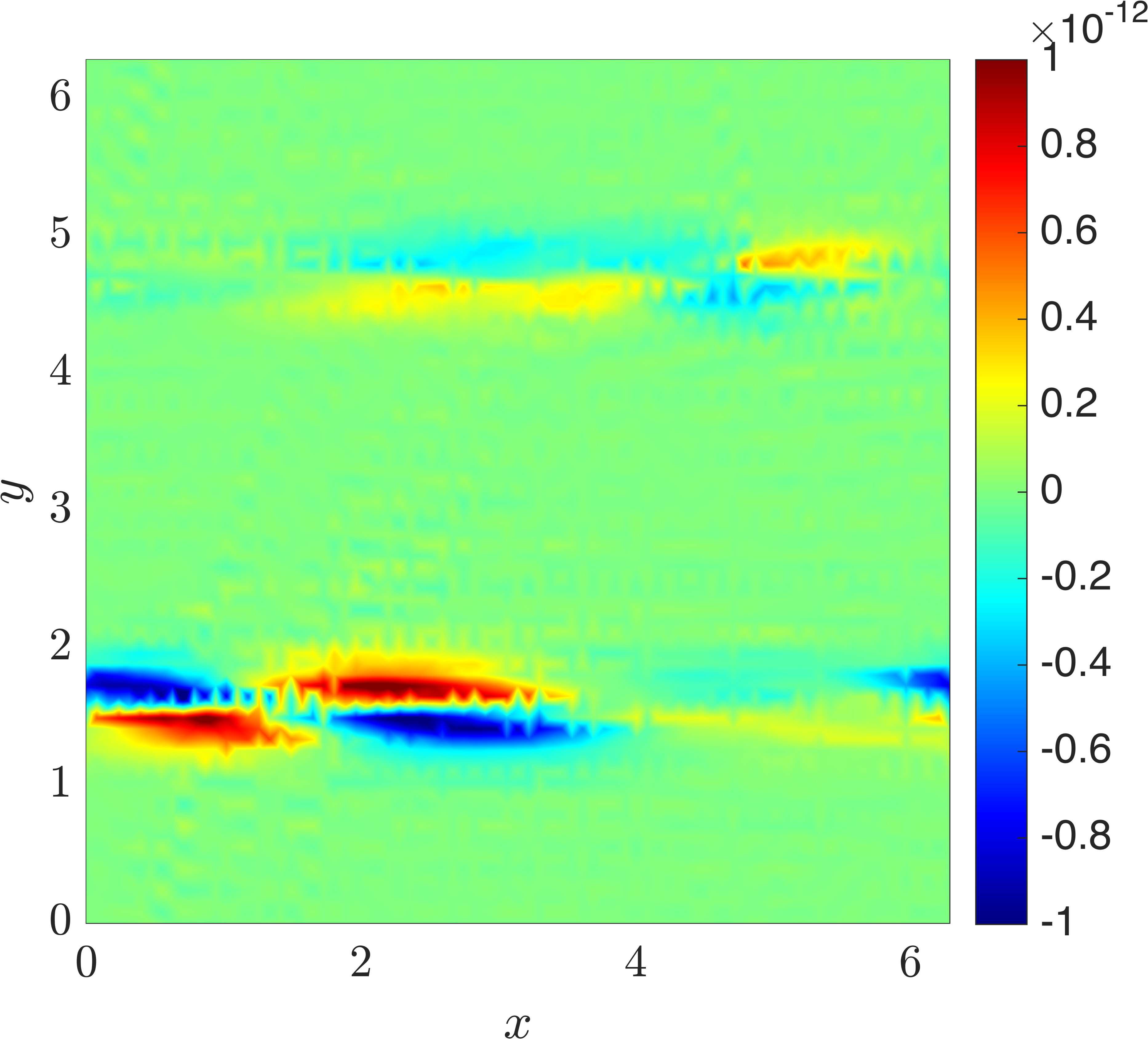}}{a) $\Delta t = \frac{1}{2}$.}
				\subsubfloat{\includegraphics[width=0.45\textwidth]{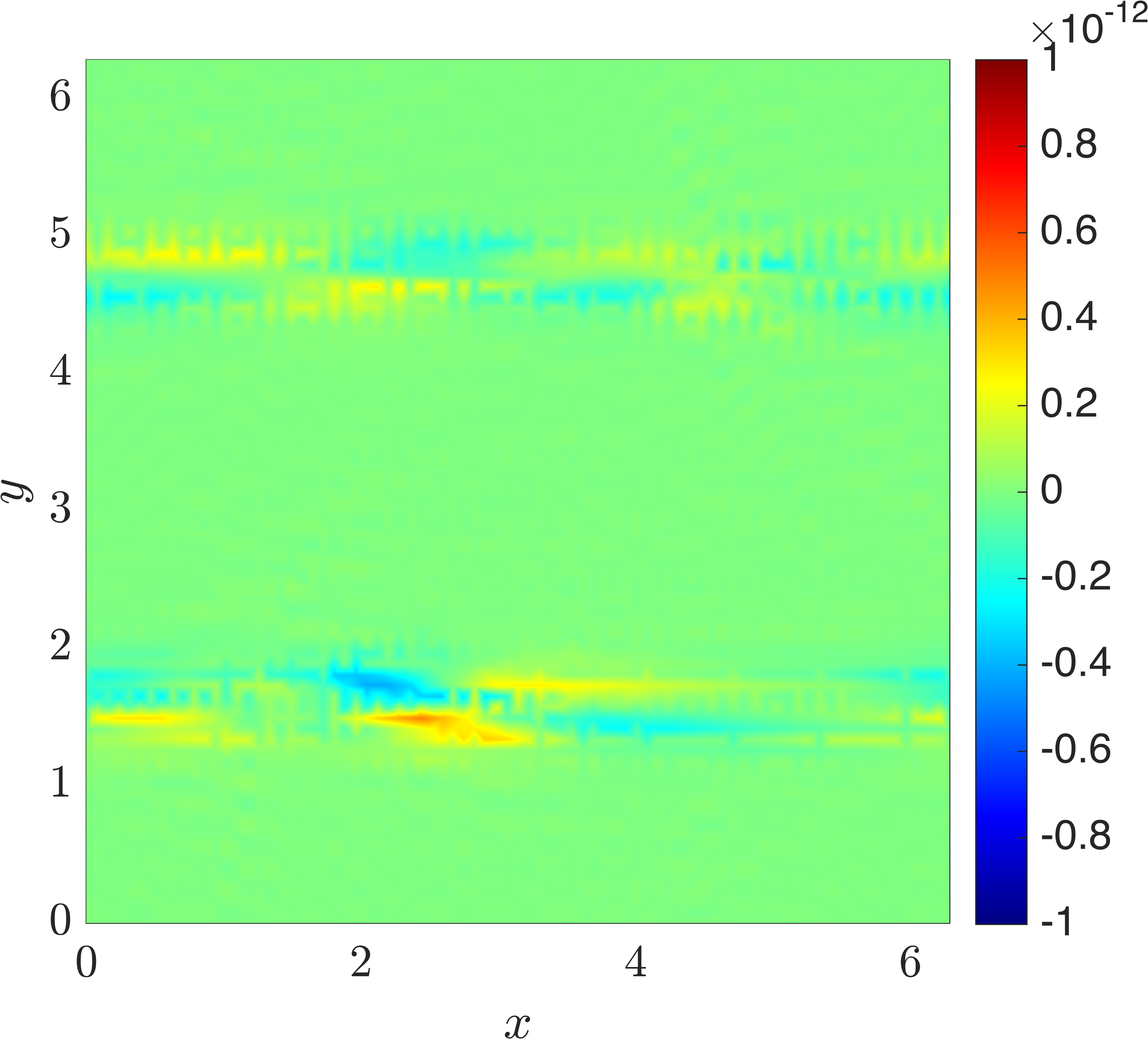}}{b) $\Delta t = 1$.}
			\end{minipage}			
			
			\caption{Error between the initial vorticity field and the vorticity field computed by evolving the inviscid Navier-Stokes equations to $t=8$s and then reversing the time step until $t=0$s is reached again. The results shown correspond to a mesh of 6400 triangular finite elements, polynomial degree $p=1$, and different time steps $\Delta t = 1,\frac{1}{2}, \frac{1}{4},\frac{1}{8}$ (starting from the bottom right and going clockwise).}
			\label{fig:time_reversibility_plots}
		\end{figure}

		The final test addresses time reversibility. In order to investigate this property of the solver, we let the flow evolve from $t=0$ to $t=8$ and then reversed the time evolution and evolved again for the same time, corresponding to the evolution from $t=8$ to $t=0$. We have performed this study with 6400 triangular finite elements of polynomial degree $p=1$ and different time steps $\Delta t = 1,\frac{1}{2}, \frac{1}{4},\frac{1}{8}$. In  \figref{fig:time_reversibility_plots} we show the error between the initial vorticity field and the final vorticity field. As can be seen, the errors are in the order of $10^{-12}$, showing the reversibility of the method up to machine accuracy.
		
	\FloatBarrier

\section{Conclusions} \label{sec::conclusions}
	This paper presents a new arbitrary order finite element solver for the Navier-Stokes equations. The advantages of this method are: (i) high order character of the spatial discretization, (ii) geometric flexibility, allowing for unstructured triangular grids, (iii) exact preservation of conservation laws for mass, kinetic energy, enstrophy and vorticity, and (iv) time reversibility.
	
	The construction of a numerical flow solver with these properties leads to a very robust method capable of performing very under-resolved simulations without the characteristic blow up of standard discretizations.

	For a simpler Taylor-Green analytical solution, the convergence properties of the proposed method have been shown. With the more challenging shear layer roll-up test case, we have shown the robustness of the method. Conservation of kinetic energy, enstrophy and vorticity up to machine accuracy was shown for this test case, even for low order approximations. 

	Nevertheless, the proposed method still allows further improvement, mainly in terms of time integration. One of the advantages of the method is the fact that the time integration method avoids the solution of a fully coupled nonlinear system of equations. Instead two quasi-linear systems of equations are solved, without the need for expensive iterations for each time step. The downside of this approach is that, for now, the time integration convergence rate is limited to first order. In the future we intend to compare this approach to a fully implicit formulation using higher order Gauss integration. Another aspect which is currently being investigated is the treatment of solid boundaries.

\bibliographystyle{elsarticle-num} 
\def\url#1{}
\bibliography{./library_clean}

\end{document}